\documentclass[A4paper]{article}

\usepackage{style}

\newcommand{\footremember}[2]{%
    \footnote{#2}
    \newcounter{#1}
    \setcounter{#1}{\value{footnote}}%
}
\newcommand{\footrecall}[1]{%
    \footnotemark[\value{#1}]%
} 

\begin{document}


\title{Improving high-order VEM stability on badly-shaped elements}

\author{S. Berrone\footremember{trailer}{Department of Applied Mathematics, Politecnico di Torino, Italy
  (stefano.berrone@polito.it, gioana.teora@polito.it, fabio.vicini@polito.it).}, G. Teora\footrecall{trailer}{}, F. Vicini\footrecall{trailer}{}}

\maketitle

\begin{abstract}
For the 2D and 3D Virtual Element Methods (VEM), a new approach to improve the conditioning of local and global matrices in the presence of badly-shaped polytopes is proposed. It defines the local projectors and the local degrees of freedom with respect to a set of scaled monomials recomputed on more well-shaped polytopes. This new approach is less computationally demanding than using the orthonormal polynomial basis.
The effectiveness of our procedure is tested on different numerical examples characterized by challenging geometries of increasing complexity.
\end{abstract}

\textbf{Keywords}: 

Ill-conditioning, Virtual Element Method, Polygonal mesh, Polyhedral mesh, Polynomial basis

\section{Introduction}

In recent years, numerical methods for the approximation of Partial Differential Equations (PDEs) on polygonal and polyhedral meshes, such as Virtual Element Methods (VEM) \cite{LBe13,LBe14} or Hybrid High Order (HHO) methods \cite{basicHHO} are gaining considerable interest since offer a convenient framework to handle challenging geometries. 
In particular, the Virtual Element Method, introduced in \cite{LBe13} for the Poisson problem and then extended in \cite{LBe141} to general second-order elliptic problems with variable coefficients, is a generalization of the Finite Element Method (FEM), which includes suitable non-polynomial functions as well as the usual polynomial functions in the local space, to employ generic polytopal meshes and to build high-order methods.
The use of these features is made possible by the introduction of suitable projection operators and by the careful selection of the local space and of local degrees of freedom (DOFs), which eliminate the need to compute in a closed form these non-polynomial functions.
It is known from VEM literature (see for example \cite{SBe17}) that the resulting system matrix is ill-conditioned in presence of badly-shaped elements (collapsing bulks, small edges,...) when resorting to the scaled monomial basis in the definition of both the local projectors and the local DOFs. 

In \cite{SBe17, mascotto, 3Dortho}, it has been suggested to replace the scaled monomial basis with an orthonormal polynomial basis to cure ill-conditioning and make the VEM solution more reliable and accurate, but this strategy can be very expensive from a computational point of view. In this paper, we propose an alternative strategy to the use of an orthonormal polynomial basis, which is much less expensive and has already led to an improvement of global performances in the two-dimensional setting of HHO \cite{inertiaHHO}.  
It consists of recomputing the scaled monomial basis on suitable polytopes whose inertia tensor is the identity tensor (re-scaled by a proper constant) and, thus of defining the local projectors and the local degrees of freedom as a function of such new polynomial basis in order to limit the condition numbers of local matrices and to improve global performances.

The structure of this work is as follows.
In Section~\ref{sec:mapping}, we define the desired properties of a well-shaped polytope and we build an affine isomorphism that allows transforming a generic polytope into a new one that has the requested features. 
In Section~\ref{sec:vemdiscretization}, we introduce the model problem and the VEM discretization on well-shaped polytopes.
Finally, in Section~\ref{sec:numericalexperiments}, we propose some numerical experiments that show the advantages of using the new procedure with respect to the standard monomial basis or the orthonormal basis both in the two and three-dimensional cases. 
The strategy presented in this work is designed and tested in the case of convex polytopes. The treatment of concave polytopes paves the way to a huge number of situations that we do not analyze.

Throughout this paper, we use the following notations. Given a polytope $E$, we denote by $h_E = \max_{\x,\bm{y} \in E} \| \x - \bm{y} \|$, $\x_E$ and $\vert E \vert$ its diameter, centroid and  measure (i.e. length or area or volume), respectively. Let $\Omega \subset \R^d$, with $d =2,3$, be a bounded polytopal domain. 
We consider a decomposition $\Th$ of $\Omega$ made of polytopal elements $E$, where we fix, as usual, $h = \max_{E\in\Th} h_E$.
We further denote by $N_v^E$, $\Eh[E]$ and $N_e^E$  the number of vertices, the set of edges and the number of edges of $E \in \Th$, respectively.  In addition, if $d=3$, $\Fh[E]$ indicates the set of the $N_f^E = \# \Fh[E]$ faces of the polyhedron $E \in \Th$ and we set $\Fh = \cup_{E \in \Th} \Fh[E]$.
Moreover, we denote by $\Poly[d]{k}{E}$ the set of polynomials defined on $E$ of degree less or equal to $k \geq 0$ and by $n^d_k = \dim \Poly[d]{k}{E}$. For the ease of notation, we further set $\Poly[d]{-1}{E} = \{0\}$ and $n^d_{-1} = 0$ and we use the two natural functions $\ell_{2}: \mathbb{N}^2 \leftrightarrow \mathbb{N}$ and $\ell_{3}: \mathbb{N}^3 \leftrightarrow \mathbb{N}$ such that:

\begin{align}
  \nonumber \left(  0, 0 \right) \leftrightarrow 1, \quad  \left(  1, 0 \right) &\leftrightarrow 2, \quad  \left(  0, 1 \right) \leftrightarrow 3, \quad \left(  2, 0 \right) \leftrightarrow 4, \dots\\ \label{eq:ell} \\
   \nonumber \left(  0, 0 ,0 \right) \leftrightarrow 1, \ \left(  1, 0 , 0 \right) \leftrightarrow 2,& \  \left(  0, 1 ,0\right) \leftrightarrow 3, \ \left(  0, 0,1 \right) \leftrightarrow 4,  \ \left(  2, 0,0 \right) \leftrightarrow 5, \dots   
\end{align}
Finally, as usual, we use $\scal[w]{}{}$ and $\norm[w]{}$ to indicate the inner product and the norm in the Lebesgue space $\lebl{w}$ on some open subset $\omega \subset \R^d$, respectively.

\section{Scaled monomials on well-shaped polytopes}
\label{sec:mapping}
We list in what follows the properties of a well-shaped polytope $E$, where we define the set of scaled monomials of degree less or equal to $k$, with $k \geq 0$, i.e. the set

\begin{equation*}
    \M[d]{k}{E} = \left\{ m^{k,d}_{\alpha} = \left(\frac{\x-\x_E}{h_E}\right)^{\bm{\alpha}}: \bm{\alpha} = \ell_d(\alpha) \in \mathbb{N}^d, \alpha = 1,\dots, n^d_k\right\},
    \label{eq:setScaledMonomial}
\end{equation*}
where $\ell_d$ is the function introduced in equation~\eqref{eq:ell}.
Let us denote by $\TT^E \in \R^{d\times d}$ the inertia tensor of $E$ with respect to its centroid $\x_E$ and the $x_i$-axes with $i=1,\dots,d$.
We further define the \textit{anisotropic ratio} of $E$ as the quotient

\begin{equation*}
    r^E = \frac{\mu_{\max}^E}{\mu_{\min}^E},
\end{equation*}
where $\mu_{\max}^E$, $\mu^E_{\min}$ are the maximum and the minimum eigenvalue of $\TT^E$, respectively.

Firstly, we recall that the products of inertia, i.e. the extra-diagonal entries of $\TT^E$, represent a measure of the imbalance in the mass distribution.
Secondly, as in \cite{anisotropic}, we say that the polytope $E$ is \emph{isotropic} if $r^E \approx 1$. 
Finally, $E$ is a \emph{well-shaped} polytope if its tensor of inertia $\TT^E$ is a diagonal matrix and $E$ is an isotropic polytope.

Thus, the goal is to define an affine isomorphism $F_E$ for each element $E\in\Th$ such that $F_E^{-1}$ maps $E$ into a well-shaped polytope $\srescale{E}$.

\subsection{2D Mapping}\label{sec:2Dmapping}

Let us consider $d=2$. Given a polygon $E \in \Th$, the inertia tensor associated to $E$ with respect to its centroid $\x_E$ and the $x_1$, $x_2$-axes is
\begin{equation}
    \TT^E = {\small \begin{bmatrix}
        \int_E (x_2-(\x_E)_2)^2               & - \int_E (x_1-(\x_E)_1)(x_2-(\x_E)_2) \\
        - \int_E (x_1-(\x_E)_1)(x_2-(\x_E)_2) & \int_E (x_1-(\x_E)_1)^2
    \end{bmatrix}}.
    \label{eq:inrtiatensor_2d}
\end{equation}
Furthermore, we define the mass matrix related to $E$ as
\begin{equation}
    \label{eq:massMatrix}
    \HH^E = \int_E (\x-\x_E)(\x-\x_E)^T \in \R^{d\times d},
\end{equation}
that is a symmetric positive-definite real matrix, and, we consider its spectral decomposition, i.e.

\begin{equation*}
    \HH^E = \QQ^E \LLambda^E \left(\QQ^E\right)^T,
\end{equation*}
where $\mathbf{Q}^E \in \R^{d\times d}$ is the orthonormal matrix whose columns represent the eigenvectors of the mass matrix $\HH^E$ and $\mathbf{\Lambda}^E$ is the diagonal matrix whose diagonal entries are the eigenvalues $\lambda_{i}^E$, $i=1,\dots,d$ of $\HH^E$.
Thus, we define a new element $\inertia{E}$ through the affine map

\begin{equation}
    \inertia{\x} \mapsto \BB^E \left(\x - \x_E\right),
    \label{eq:inertiamap}
\end{equation}
where $\BB^E = \sqrt{\lambda_{\max}^E} \sqrt{\mathbf{\left(\Lambda^E\right)}^{-1}} \left(\mathbf{Q}^E\right)^T$ is invertible and such that $\abs{\det \BB^E} = \frac{\left(\lambda_{\max}^E\right)^{d/2}}{\sqrt{\prod_{i=1}^d \lambda_i^E}}$,  $\lambda_{\max}^E = \max_{i=1,\dots,d} \lambda_i^E$.
We note that the centroid of $\inertia{E}$ is $\inertia{\x}_{\inertia{E}} = \bm{0}$ and $\TT^{\inertia{E}}$ is a diagonal matrix. Indeed,

\begin{align*}
    \TT^{\inertia{E}}_{ij} & = - \int_{\inertia{E}} \inertia{x}_i \inertia{x}_j                                                                                                                         \\
                           & = - \abs{\det \BB^E}  \BB^E(i,:) \int_E (\x - \x_E)(\x-\x_E)^T \left(\BB^E(j,:)\right)^T                                                                                                  \\
                           & = - \abs{\det \BB^E} \BB^E(i,:) \HH^E\ \left(\BB^E(j,:)\right)^T                                                                                                           \\
                           & = - \lambda_{\max}^E \abs{\det \BB^E} \frac{\left(\QQ^E(:,i)\right)^T}{\sqrt{\lambda_i^E}} \QQ^E \LLambda^E \left(\QQ^E\right)^T \frac{\QQ^E(:,j)}{\sqrt{\lambda_j^E}} \\
                           & = - \lambda_{\max}^E \abs{\det \BB^E} \left(\frac{\left(\QQ^E(:,i)\right)^T}{\sqrt{\lambda_i^E}} \left(\QQ^E \sqrt{\LLambda^E}\right) \right) \left(\left(\QQ^E \sqrt{\LLambda^E}\right)^T \frac{\QQ^E(:,j)}{\sqrt{\lambda_j^E}}\right) \\
                           & = - \lambda_{\max}^E \abs{\det \BB^E} \ee_i^T \ee_j = 0,\quad \forall i,j =1,\dots,d \text{ s.t. } i \neq j,
\end{align*}
where, given a generic matrix $\mathbf{A}$, $\mathbf{A}(i,:)$ is the sub-matrix of $\mathbf{A}$ made up of its $i$-th row, $\mathbf{A}(:,j)$ is the sub-matrix of $\mathbf{A}$ made up of its $j$-th column and the set $\{\ee_i\}_{i=1}^d$ represents the canonical basis of $\R^d$.
Finally, concerning the diagonal entries of $\TT^{\inertia{E}}$, we note that

\begin{align*}
    \TT^{\inertia{E}}_{ss} = \int_{\inertia{E}} \inertia{x}_i \inertia{x}_i
     & = \abs{\det \BB^E} \BB^E(i,:) \HH^E\ \left(\BB^E(i,:)\right)^T \\
     & = \lambda_{\max}^E \abs{\det \BB^E} \ee_i^T \ee_i              \\
     & = \lambda_{\max}^E \abs{\det \BB^E},\quad \forall s,i=1,2,\ i \neq s.
\end{align*}
Thus, the diagonal elements prove to be constant with respect to the matrix index $s$. In conclusion, the new element $\inertia{E}$ is isotropic and well-shaped according to the definitions we provide.

The computational cost of this mapping depends only on the dimension $d$ of the problem but not on the order of accuracy of the discretization method. For stability reasons, in order to avoid small eigenvalues, it is preferable to perform a re-scaling before proceeding with the mapping \eqref{eq:inertiamap}. Furthermore, we also decide to re-scale elements after the application of the mapping \eqref{eq:inertiamap} in order to have polygons with unit diameter.

In conclusion, on each element $E\in\Th$, we perform sequentially the transformations

\begin{equation}
\centering
\begin{aligned}
    E \longrightarrow &\frescale{E} \longrightarrow \inertia{E} \longrightarrow \srescale{E} \\
    \x \mapsto \frescale{\x} = \frac{1}{h_{E}} \x \mapsto \inertia{\x} &= \BB^{\frescale{E}} \left(\frescale{\x} - \frescale{\x}_{\frescale{E}}\right) \mapsto \srescale{\x} =\frac{1}{h_{\inertia{E}}} \inertia{\x}. 
\end{aligned}
\label{eq:mapping}
\end{equation}
and we define the local mapping $F_E$ such that
\begin{equation}
    \label{eq:mapFe}
    \x = F_E(\srescale{\x}) = \x_E + \FF^E\srescale{\x},
\end{equation}
where

\begin{equation*}
    \FF^E =  {h_{E}}  {h_{\inertia{E}}} \left(\BB^{\frescale{E}}\right)^{-1} \text{ and } \abs{ \det \FF^E} = h_{E}^d h_{\inertia{E}}^d  \abs{ \det \left(\BB^{\frescale{E}}\right)^{-1} }.
\end{equation*}
Examples of $F_E$ mapping are shown in Table~\ref{tab:Fmap_2D}.

\begin{table}[ht]
    \resizebox{\textwidth}{!}{
        \begin{tabular}
            {cc|cc}
            {\includegraphics[width=1.5in]{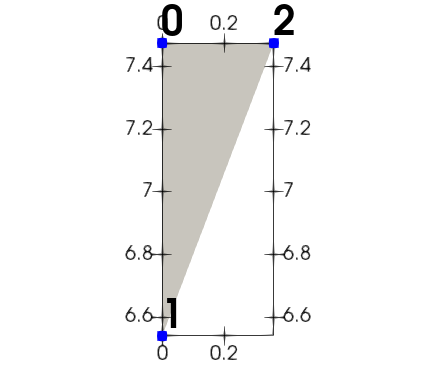}}          &
            {\includegraphics[width=1.5in]{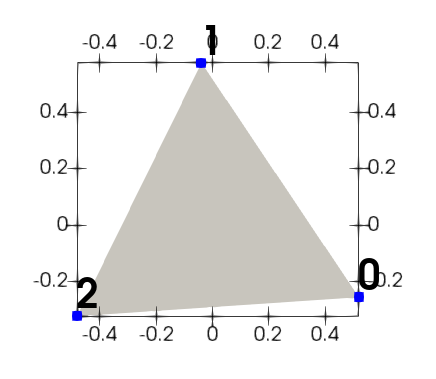}}  &
            {\includegraphics[width=1.5in]{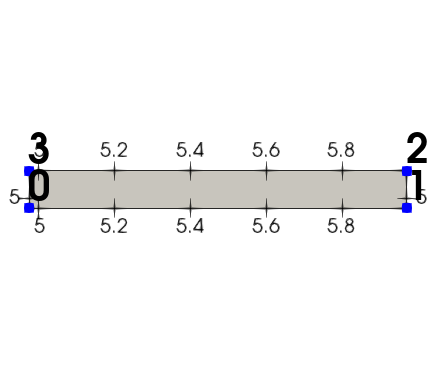}}          &
            {\includegraphics[width=1.5in]{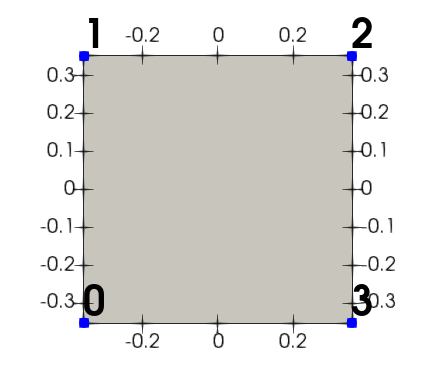}}                                          \\
            $E$                                                                     & $\srescale{E}$ & $E$ & $\srescale{E}$ \\ \midrule
            {\includegraphics[width=1.5in]{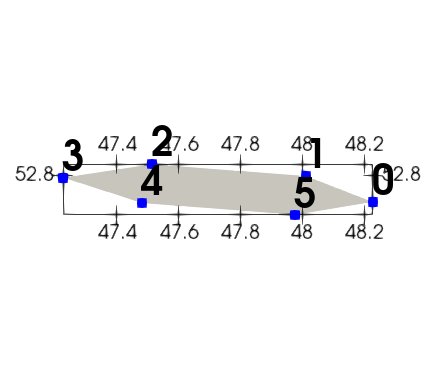}}         &
            {\includegraphics[width=1.5in]{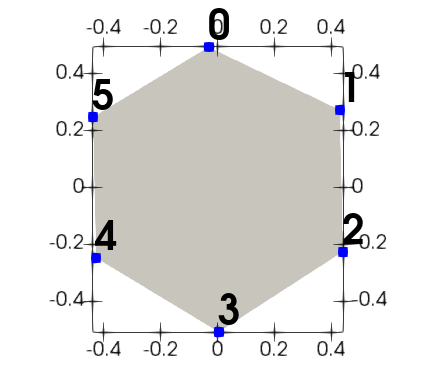}} &
            {\includegraphics[width=1.5in]{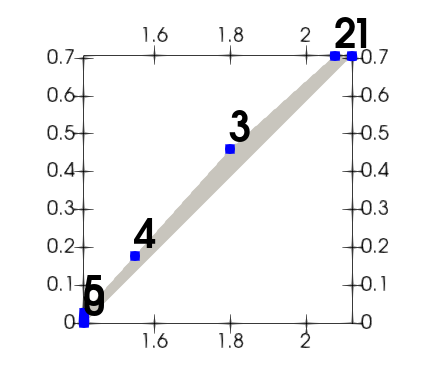}}           &
            {\includegraphics[width=1.5in]{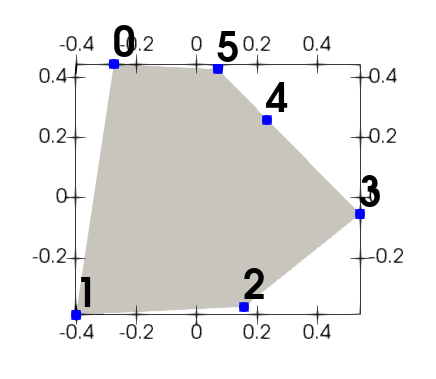}}                                           \\
            $E$                                                                     & $\srescale{E}$ & $E$ & $\srescale{E}$
        \end{tabular}
    }
    \caption{Examples of $F_E$ mapping if $d=2$.}
    \label{tab:Fmap_2D}
\end{table}

\subsection{3D Mapping}

Let us set $d=3$.
Given a polyhedron $E \in \Th$, we define the inertia tensor $\TT^E$ related to $E$ with respect to its centroid $\x_E$ and $x_1$, $x_2$, $x_3$-axes as the tensor whose entries are defined as

\begin{equation*}
\TT^{E}_{ss}  = \int_{E} (x_i-(\x_E)_i)^2 + (x_j-(\x_E)_j)^2 ,  \quad \forall s,i,j =1,\dots,d,\ s\neq i\neq j,
\end{equation*}
and

\begin{equation*}
\TT^{E}_{ij}  = -\int_{E} (x_i-(\x_E)_i)(x_j-(\x_E)_j) ,  \quad \forall i,j =1,\dots,d,\ i\neq j.
\end{equation*}
A natural extension to the 3D case of the affine map \eqref{eq:mapFe} is here defined by exploiting the spectral decomposition of the tridimensional mass matrix $\HH^E$.

We first note that the resulting polyhedron $\inertia{E}$ is a well-shaped polyhedron. 
Indeed,

\begin{align*}
    \TT^{\inertia{E}}_{ss} & = \int_{\inertia{E}} \inertia{x}_i \inertia{x}_i + \int_{\inertia{E}} \inertia{x}_j \inertia{x}_j                          \\
               & = \abs{\det \BB^E} \left[\BB^E(i,:) \HH^E\ \left(\BB^E(i,:)\right)^T + \BB^E(j,:) \HH^E\ \left(\BB^E(j,:)\right)^T \right] \\
               & = \lambda^E_{\max} \abs{\det \BB^E}  \left[\ee_i^T \ee_i + \ee_j^T \ee_j \right]                                           \\
               & = 2\lambda^E_{\max} \abs{\det \BB^E} ,  \quad \forall s,i,j =1,2,3,\ s\neq i\neq j.
\end{align*}
Moreover, as highlighted in \cite{LBe14}, if $d=3$, we need to compute both the 3D projectors on $E$ and the 2D projectors on each face $f \in \Fh[E]$, thus we need to define a polynomial basis $\Poly[3]{k}{E}$ $\forall E \in \Th$ and a polynomial basis $\Poly[2]{k}{f}$ $\forall f \in \Fh$. 
For this reason, we built three different approaches (\ref{itm:B}, \ref{itm:F} and \ref{itm:BF}) by defining, $\forall E \in \Th$,
\begin{description}[style=multiline, labelwidth=1cm]
    \item[\namedlabel{itm:B}{(B)}] the 3D set of scaled monomials on the well-shaped polyhedron $\srescale{E} = F_E^{-1}(E)$ and the 2D sets of scaled monomials on the original faces $f \in \Fh[E]$;
    \item[\namedlabel{itm:F}{(F)}] the 3D set of scaled monomials on the original polyhedron $E$ and the 2D sets of scaled monomials on the mapped faces $\mf{f}$ obtained by applying a mapping $F_f^{-1}$ to each face $f\in \Fh[E]$ defined accordingly to what is done for the bidimensional case in the Section \ref{sec:2Dmapping};
    \item[\namedlabel{itm:BF}{(B-F)}] the 3D set of scaled monomials on $\srescale{E} = F_E^{-1}(E) $ and the 2D sets of scaled monomials on the mapped faces $\mf{f} = F_f^{-1}(f)$ for each $f\in \Fh[E]$.
\end{description}
We want to highlight that in the \ref{itm:BF} approach the new polygons $\mf{f}$ are obtained by applying $F_f^{-1}$ to the original faces $f$ of $E$ and not to faces $\srescale{f}$ of $\srescale{E}$.

Examples of the maps $F_E$ are shown in Table \ref{tab:Fmap_3D}.

\begin{table}[ht]
    \centering
    \begin{tabular}{cc}
        {\includegraphics[width=1.5in]{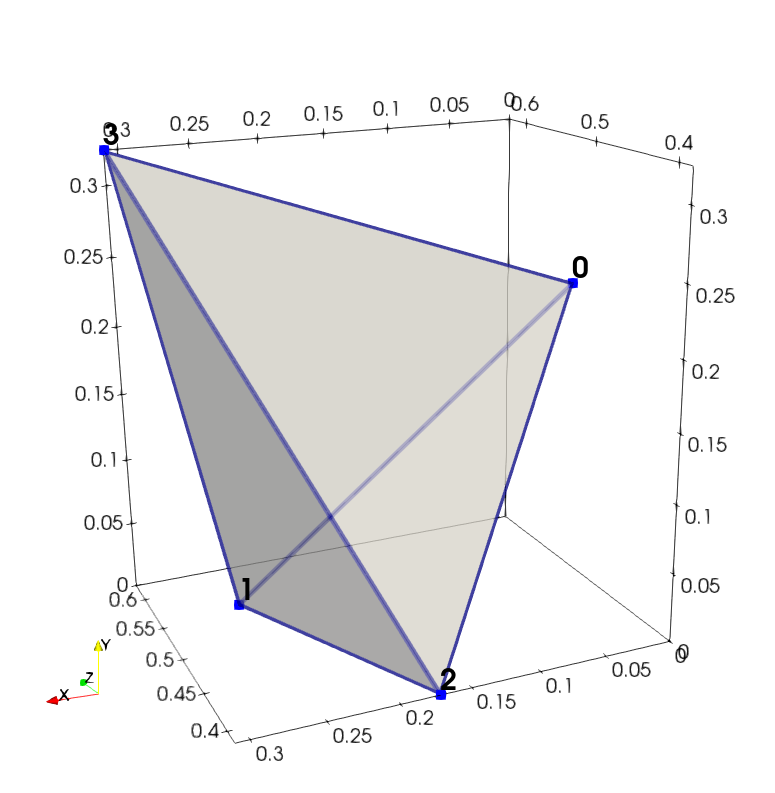}}   &
        {\includegraphics[width=1.5in]{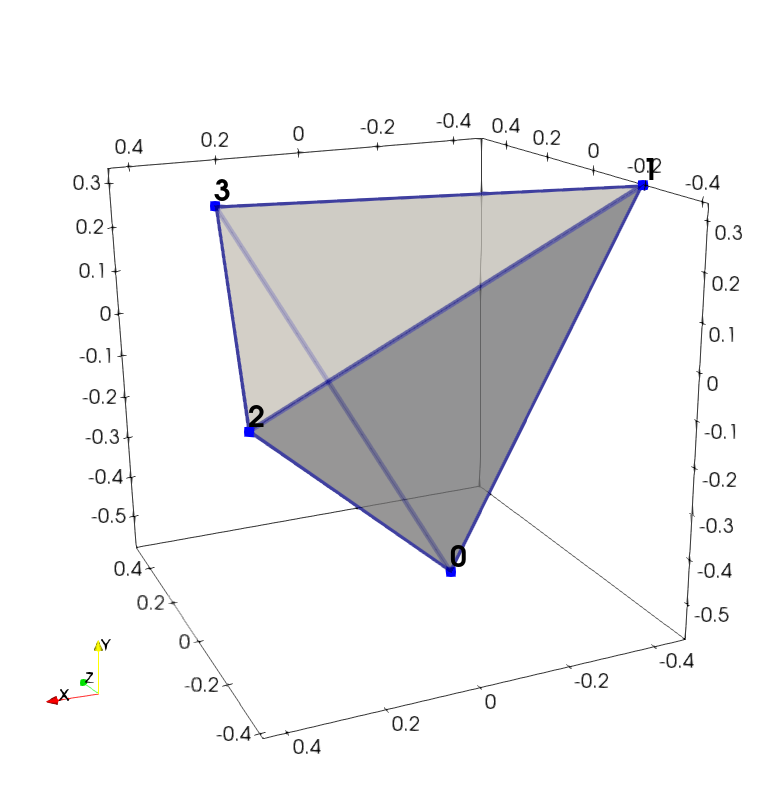}}            \\
        $E$                                                            & $\srescale{E}$ \\ \hline
        {\includegraphics[width=1.5in]{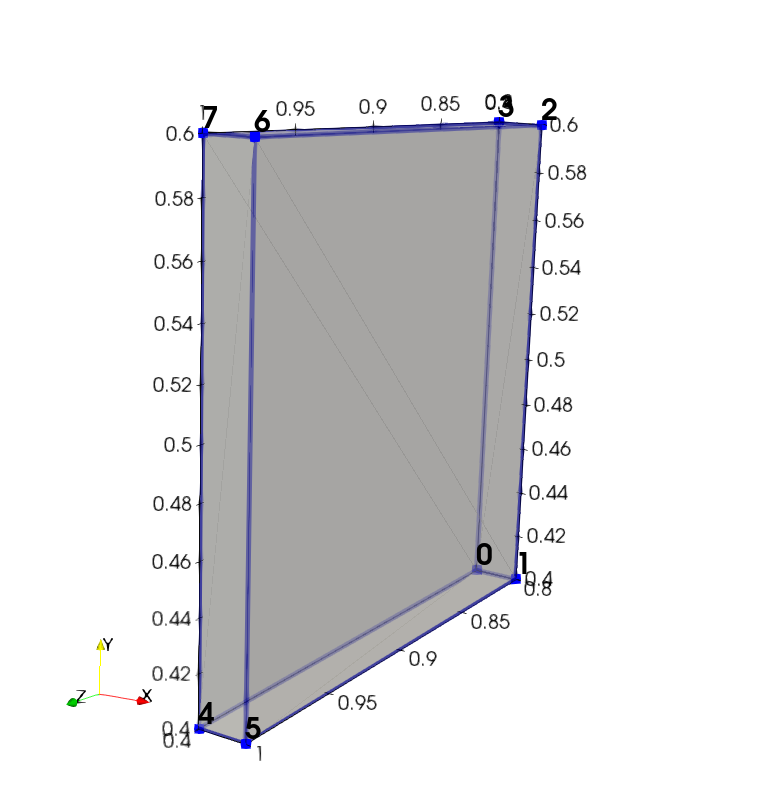}} &
        {\includegraphics[width=1.5in]{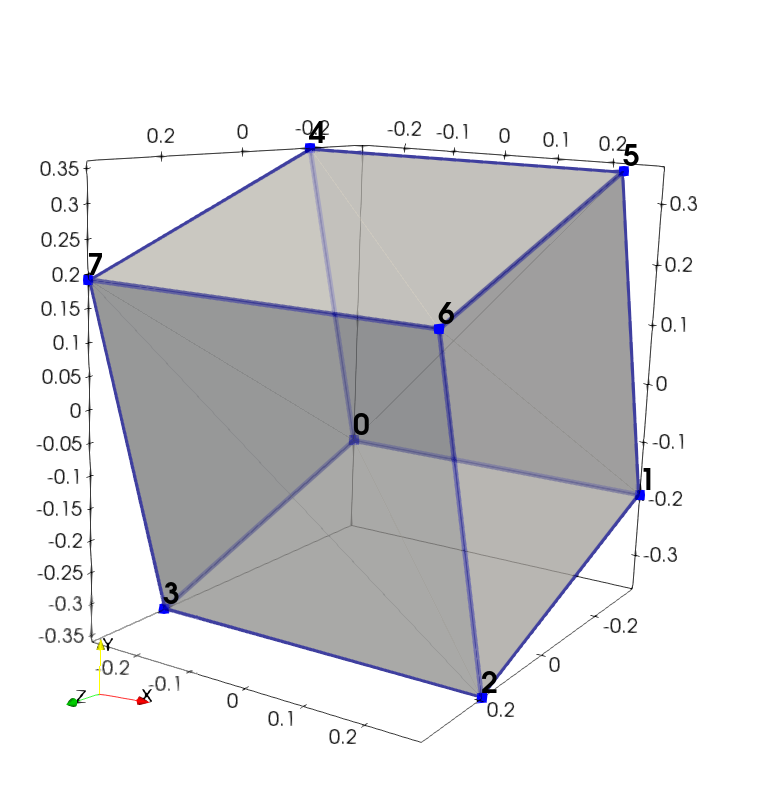}}          \\
        $E$                                                            & $\srescale{E}$
    \end{tabular}
    \caption{Examples of $F_E$ mapping if $d=3$.}
    \label{tab:Fmap_3D}
\end{table}
\section{VEM discretization}\label{sec:vemdiscretization}

In the following, for the sake of convenience, if $d=3$, we define $F_f$, $\forall f \in \Fh$, in the \ref{itm:B} approach and $F_E$, $\forall E\in\Th$, in the \ref{itm:F} approach as the identity maps.

Given $k \geq 1$, on each element $\srescale{E} = F_E^{-1}(E)$ with $E\in\Th$, we introduce the local projection operators $\mproj{\srescale{\nabla},\srescale{E}}{k}: H^1(\srescale{E})   \rightarrow \Poly[d]{k}{\srescale{E}}$ and $\mproj{0,\srescale{E}}{k}: \lebl{\srescale{E}} \rightarrow \Poly[d]{k}{\srescale{E}}$, as
\begin{equation}
    \begin{cases} \scal[\srescale{E}]{ \srescale{\nabla} \srescale{p}}{\srescale{\nabla}  \left(  \mproj{\srescale{\nabla},\srescale{E}}{k} \srescale{v} - \srescale{v} \right)} =0, & \forall \srescale{p} \in \Poly[d]{k}{\srescale{E}} \\
        \srescale{P}_{0} \left(  \mproj{\srescale{\nabla},\srescale{E}}{k} \srescale{v} - \srescale{v}  \right) =0
    \end{cases},
    \label{eq:PiNablak_cond_ortho}
\end{equation}
and
\begin{equation}
    \scal[\srescale{E}]{ \srescale{p}}{ \left( \mproj{0,\srescale{E}}{k} \srescale{v} - \srescale{v} \right)  }=0,\quad \forall \srescale{p} \in \Poly[d]{k}{\srescale{E}},
    \label{eq:Pi0km1_cond_ortho}
\end{equation}
where
\begin{equation}
    \srescale{P}_0 (\srescale{v}) = \begin{cases}
        \left(\srescale{v},1\right)_{\partial \srescale{E}} & \text{if } k=1,\\
        \left(\srescale{v},1\right)_{ \srescale{E}}         & \text{if } k > 2.
    \end{cases}
\end{equation}
Furthermore, for the sake of convenience, we use the symbol $\mproj{0,\srescale{E}}{k-1}$ also for the $L^2$-projection operator of vector-valued functions onto the polynomial space $\left(\mathbb{P}_{k-1}(\srescale{E})\right)^2$, meaning a component-wise application.

Now, on each $E \in \Th$, following \cite{LBe141,BAh13}, we introduce the local enhanced Virtual Element spaces:
\begin{itemize}
    \item if $d=2$,
          \begin{align*}
              V_k^2(E)\! = & \big\{ v\!\in\! H^1(E)\!: \Delta v\!\in\! \mathbb{P}_k(E), v_{|e}\!\in\!\mathbb{P}_k(E)\ \forall e\!\in\! \mathcal{E}_{h,E}, v_{\vert \partial E} \!\in\! C^0(\partial E),\\
                                    & \scal[\srescale{E}]{\srescale{v}}{\srescale{p}} = \scal[\srescale{E}]{\mproj{\srescale{\nabla},\srescale{E}}{k}\srescale{v}}{\srescale{p}}\ \forall \srescale{p} \in \Poly[2]{k}{\srescale{E}}/\ \Poly[2]{k-2}{\srescale{E}} \big\};
          \end{align*}
    \item if $d=3$,
          \begin{align*}
              V_k^3(E)\! = & \big\{ v\!\in\! H^1(E)\!: \Delta v\!\in\! \mathbb{P}_k(E), v_{\vert f}\!\in\! V^2_k(f)\ \forall f\!\in\!\Fh[E], v_{\vert \partial E}\!\in\! C^0(\partial E),\\
                                    & \scal[\srescale{E}]{\srescale{v}}{\srescale{p}} = \scal[\srescale{E}]{\mproj{\srescale{\nabla},\srescale{E}}{k}\srescale{v}}{\srescale{p}}\ \forall \srescale{p} \in \Poly[3]{k}{\srescale{E}}/\ \Poly[3]{k-2}{\srescale{E}}  \big\},
          \end{align*}
\end{itemize}
where $\srescale{w}(\srescale{\x}) = w(\x_E + \FF^E \srescale{\x})$ and $\nabla w = \left(\FF^E\right)^{-T} \srescale{\nabla} \srescale{w}$, $\forall w \in V_k^d\left(E\right)$. We further define the following set of local DOFs: $\forall v_h \in V^d_k(\srescale{E})$
\begin{enumerate}
    \item the value of $v_h$ at the vertices of $E$;
    \item if $k > 1$, for each edge $e\in\Eh[E]$, the value of $v_h$ at the $k-1$ internal Gauss-Lobatto quadrature nodes on $e$;
    \item if $k > 1$ and $d=3$, for each face $f\in\Fh[E]$, the scaled moments on $\mf{f} = F^{-1}_f(f)$
          \begin{equation}
              \frac{1}{\vert \mf{f} \vert} \scal[\mf{f}]{\mf{v}_h}{ \mf{m}^{k-2,d-1}_{\alpha}},\quad \forall \mf{m}^{k-2,d-1}_{\alpha} \in \M[d-1]{k-2}{\mf{f}};
          \end{equation}
    \item if $k > 1$, the scaled moments on $\srescale{E} = F^{-1}_E(E)$
          \begin{equation}
              \frac{1}{\vert\srescale{E} \vert} \scal[\srescale{E}]{\srescale{v}_h}{ \srescale{m}^{k-2,d}_{\alpha} },\quad \forall \srescale{m}^{k-2,d}_{\alpha} \in \M[d]{k-2}{\srescale{E}}.
          \end{equation}
\end{enumerate}

Let it be $\Ndof[E]= \dim V^d_k\left(E\right)$, for each element $E\in\Th$, we denote by $\dof[E]{i}{}$ the operator that associates its $i$-th local degree of freedom to each sufficiently smooth function $\varphi$ and by $\{\varphi_i\}_{i=1}^{\Ndof[E]}$ the set of local Lagrangian VEM basis functions related to the defined DOFs. Furthermore, we introduce the local projection matrices  $\mmproj{\srescale{\nabla},\srescale{E}}{k} \in \R^{n^d_{k}\times \Ndof[E]}$, $\mmproj{0,\srescale{E}}{k-1} \in \R^{n^d_{k-1}\times \Ndof[E]}$ and $\mmproj{0,\srescale{x}_{j},\srescale{E}}{k-1} \in \R^{n^d_{k-1}\times \Ndof[E]}$, for $j=1,\dots,d$, which are defined as
\begin{equation}
    \!
    \begin{aligned}
        \mproj{\srescale{\nabla},\srescale{E}}{k} \varphi_i = \sum_{\alpha=1}^{n^d_{k}} \left(\mmproj{\srescale{\nabla},\srescale{E}}{k}\right)_{\alpha i} \srescale{m}^{k,d}_{\alpha},\quad \mproj{0,\srescale{E}}{k-1} \varphi_i = \sum_{\alpha=1}^{n^d_{k-1}} \left(\mmproj{0,\srescale{E}}{k-1}\right)_{\alpha i} \srescale{m}^{k-1,d}_{\alpha}, \\
        \mproj{0,\srescale{E}}{k-1} \frac{\partial \varphi_i}{\partial \srescale{x}_j} = \sum_{\alpha=1}^{n^d_{k-1}} \left(\mmproj{0,\srescale{x}_j,\srescale{E}}{k-1}\right)_{\alpha i} \srescale{m}^{k-1,d}_{\alpha},\ \forall j=1,\dots,d.
    \end{aligned}
    \label{eq:localmatrices}
\end{equation}
We note that the definition of local DOFs makes the computation of the local projection matrices \eqref{eq:localmatrices} completely independent of the geometric properties of the original element $E$ in the bidimensional case. On the other hand, in the tridimensional case, the geometric properties of the original element $E$ influence the computation of the projectors with an intensity that depends on the chosen approach (\ref{itm:B}, \ref{itm:F} or \ref{itm:BF}).

Finally, we introduce the global Virtual Element space

\begin{equation*}
    V_{h,k} = \left\{v \in C^0\left(\overline{\Omega}\right) \cap H^1_0\left(\Omega\right): \ v\in V_k^d\left(E\right)\ \forall E \in \Th \right\}.
\end{equation*}
and, in agreement with the local choice of the DOFs, we define the following set of DOFs: $\forall v \in V_{h,k}$
\begin{enumerate}
    \item the value of $v_h$ at the internal vertices of the decomposition $\Th$;
    \item if $k > 1$, for each internal edge $e\in \Eh$, the value of $v_h$ at the $k-1$ internal Gauss-Lobatto quadrature nodes on $e$;
    \item if $k > 1$ and $d=3$, for each internal face $f\in \Fh$, the scaled moments on $\mf{f} = F^{-1}_f(f)$
          \begin{equation}
              \frac{1}{\vert \mf{f} \vert} \scal[\mf{f}]{\mf{v}_h}{ \mf{m}^{k-2,d-1}_{\alpha}},\quad \forall \mf{m}^{k-2,d-1}_{\alpha} \in \M[d-1]{k-2}{\mf{f}};
              \label{eq:facemoments}
          \end{equation}
    \item if $k > 1$, for each element $E\in\Th$, the scaled moments on $\srescale{E} = F^{-1}_E(E)$
          \begin{equation}
              \frac{1}{\vert\srescale{E} \vert} \scal[\srescale{E}]{\srescale{v}_h}{ \srescale{m}^{k-2,d}_{\alpha} },\quad \forall \srescale{m}^{k-2,d}_{\alpha} \in \M[d]{k-2}{\srescale{E}}.
          \end{equation}
\end{enumerate}

\begin{remark}\label{rem:notmappedfaces}
    It is very important in the approach \ref{itm:B} to set out monomials on the original faces and not on the faces of the mapped element $F_E^{-1}(E)$ in order to define uniquely the degrees of freedom \eqref{eq:facemoments} on each face $f\in \Fh$. For the same reason, in the approach \ref{itm:BF}, we apply the map $F_f$ to the original faces $f$ and not to faces $\srescale{f}$ belonging to the mapped element $\srescale{E}$.
\end{remark}

\subsection{Example: an advection-diffusion-reaction problem}
\label{sec:adrProblem}

Let $\D$ be a symmetric uniformly positive definite tensor over $\Omega$, $\gamma$ be a sufficiently smooth function $\Omega \to \R$ and $\bb$ be a smooth vector-valued function $\Omega \to \R^d$ s.t. $\nabla \cdot \bb = 0$. Given $f\in\lebl{\Omega}$, we consider the following advection-diffusion-reaction problem
\begin{equation}
    \begin{cases}
        - \nabla \cdot \left(\D\nabla u\right) + \bb \cdot \nabla u + \gamma u = f & \text{in } \Omega          \\
        u=0                                                                        & \text{on } \partial \Omega
    \end{cases}.
    \label{eq:modelproblem}
\end{equation}
Without loss of generality, we assume homogeneous boundary Dirichlet condition. The non-homogeneous case can be treated with the standard lifting procedure.

We introduce the local bilinear form, $\forall E \in \Th$

\begin{equation*}
    a^E(u,v) = \int_{E} \D \nabla u \cdot \nabla v + \int_{E} \bb \cdot \nabla u v + \int_{E} \gamma u v
\end{equation*}
and we write the variational formulation of \eqref{eq:modelproblem} as: \textit{Find} $u \in V = H^1_0(\Omega)$ \textit{such that}
\begin{equation}
    \sum_{E\in \Th} a^E(u,v) = \sum_{E\in \Th}  \int_{E} fv, \quad \forall v \in V.
    \label{eq:PV}
\end{equation}
We note that

\begin{equation*}
    \int_E \left(\D \nabla u\right) \cdot \nabla v  = \int_{\srescale{E}} \left(\srescale{\nabla} \srescale{u}\right)^T\left(\FF^E\right)^{-1}\srescale{\D} \left(\FF^E\right)^{-T} \srescale{\nabla} \srescale{v} \abs{ \det \FF^E },
\end{equation*}
\begin{equation*}
    \int_E \bb \cdot \nabla u v  = \int_{\srescale{E}} \srescale{\bb}^T\left(\FF^E\right)^{-T} \srescale{\nabla} \srescale{u} \srescale{v} \abs{ \det \FF^E },
\end{equation*}
\begin{equation*}
    \int_E \gamma u v = \int_{\srescale{E}} \srescale{\gamma} \srescale{u} \srescale{v} \abs{ \det \FF^E },\quad \int_E fv = \int_{\srescale{E}} \srescale{f} \srescale{v} \abs{ \det \FF^E }.
\end{equation*}
where, for all $\srescale{\x} \in \srescale{E}$,

\begin{equation*}
    \srescale{\D}(\srescale{\x}) = \D(\x_E + \FF^E \srescale{\x}),\quad \srescale{\bb}(\srescale{\x}) = \bb(\x_E + \FF^E \srescale{\x}),
\end{equation*}
\begin{equation*}
    \srescale{\gamma}(\srescale{\x}) = \gamma(\x_E + \FF ^E\srescale{\x}),\quad \srescale{f}(\srescale{\x}) = f(\x_E + \FF^E\srescale{\x}).
\end{equation*}

Let us introduce on each element $E$ the symmetric uniformly positive-definite tensor $\srescale{\KK} = \abs{ \det \FF^E } \left(\FF^E\right)^{-1}\srescale{\D} \left(\FF^E\right)^{-T}$, the vector-valued function $\srescale{\bbeta} = \abs{ \det \FF^E } \left(\FF^E\right)^{-1} \srescale{\bb}$ and $\srescale{\rho} = \abs{ \det \FF^E }  \srescale{\gamma}$. We define the local virtual bilinear form as: $\forall u_h, v_h \in V_{h,k}$,

\begin{align*}
    a^E_h(u_h,v_h) & = \int_{\srescale{E}} \left(\srescale{\KK} \mproj{0,\srescale{E}}{k-1}\srescale{\nabla} \srescale{u}_h\right) \cdot \mproj{0,\srescale{E}}{k-1} \srescale{\nabla} \srescale{v}_h \\&+
    C_{\D} \vemstab[E]{(I-\mproj{\srescale{\nabla},\srescale{E}}{k})\srescale{u}_h}{(I-\mproj{\srescale{\nabla},\srescale{E}}{k})\srescale{v}_h}                                                      \\
                   & + \int_{\srescale{E}} \srescale{\bbeta} \cdot \mproj{0,\srescale{E}}{k-1} \srescale{\nabla} \srescale{u}_h \ \srescale{\Pi}^{0}_{k-1} \srescale{v}_h                             \\
                   & + \int_{\srescale{E}}  \srescale{\rho}\ \mproj{0,\srescale{E}}{k-1}\srescale{u}_h \ \mproj{0,\srescale{E}}{k-1} \srescale{v}_h,
\end{align*}
where $C_{\D}$ is a constant depending on $\D$ and we employ the standard \textit{dofi-dofi} stabilization
\begin{equation}
    \vemstab[E]{u_h}{v_h} = h_E^{d-2} \sum_{i = 1}^{N^{\dof{}{}}_E} \dof[E]{i}{\srescale{u}_h} \dof[E]{i}{\srescale{v}_h}.
    \label{eq:stabterm}
\end{equation}
We note that it is very important to use the diameter $h_E$ of the original element in the equation \eqref{eq:stabterm} in order to scale $\vemstab[E]{}{}$ as $a^E(\cdot,\cdot)$. See \cite{LBe13} for further details.

Finally, the VEM discrete counterpart of \eqref{eq:PV} reads: \textit{Find} $u_h \in V_{h,k}$ \textit{such that}
\begin{equation}
    \label{eq:discreteProblem}
    \sum_{E\in \Th} a^E_h(u_h,v_h) = \sum_{E\in \Th}  \int_{\srescale{E}} \srescale{f}\ \mproj{0,\srescale{E}}{k-1}\srescale{v}_h \abs{ \det \FF^E }, \quad \forall v_h \in V_{h,k}.
\end{equation}
\section{Numerical experiments}\label{sec:numericalexperiments}

In this section, we propose some numerical experiments to validate the aforementioned approaches, which we generically called \textit{inertial} (\textit{Inrt} in short), in the two and three-dimensional cases. To show the advantages of using our procedures, we compare their performances in different discretizations of increasing complexity with respect to
\begin{itemize}
	\item the standard monomial approach (\textit{Mon} in short) described in \cite{LBe14};
	\item the orthonormal approach (\textit{Ortho} in short) which follows the construction presented in \cite{mascotto} and \cite{3Dortho} for the two and three-dimensional cases, respectively. We highlight that for the three-dimensional case, we choose to use an orthonormal basis both on the bulks and on the faces of the decomposition.
\end{itemize}
In the following, for the ease of notation, we assume to work on a mapped element $\srescale{E}$ also in the monomial and in the orthonormal approaches, where the used maps $F_{E}$, $\forall E\in \Th$, and $F_{f}$, $\forall f \in \Fh$, are the  identity maps.

The comparison is based on the analysis of the condition numbers of local projection matrices that are defined in \eqref{eq:localmatrices} and that are used to assemble the local system matrix. For the 3D case, we also analyze the condition numbers of the 2D local projection matrices $\mfmproj{\mf{\nabla},\mf{f}}{k} \in \R^{n^d_{k}\times \Ndof[f]}$, $\mfmproj{0,\mf{f}}{k-1} \in \R^{n^d_{k-1}\times \Ndof[f]}$  that are employed for computing the boundary integrals appearing in the computation of the 3D local projection matrices.
Furthermore, we analyze the behaviours of the condition number of the global system matrix $\mathbf{A}$ of the discrete problem \eqref{eq:discreteProblem} and of the following relative error norms:
\begin{equation}\footnotesize
	u_{\text{err}}^2 = {\frac{\displaystyle\sum_{E\in \Th}  \int_{\srescale{E}} \left(\srescale{u} - \mproj{0,\srescale{E}}{k} \srescale{u}_h \right)^2 \abs{ \det \FF^E }}{\displaystyle\sum_{E\in \Th} \norm[E]{ u}^2}},
	\label{eq:errorL2}
\end{equation}
\begin{equation}\footnotesize
	\nabla u_{\text{err}}^2 = {\frac{\displaystyle\sum_{E\in \Th} \int_{\srescale{E}} \left(\srescale{\nabla} \srescale{u} - \mproj{0,\srescale{E}}{k-1}\srescale{\nabla} \srescale{u}_h\right)^T \left(\FF^E\right)^{-1}\left(\FF^E\right)^{-T} \left(\srescale{\nabla} \srescale{u} - \mproj{0,\srescale{E}}{k-1}\srescale{\nabla} \srescale{u}_h\right) \abs{ \det \FF^E } }{\displaystyle\sum_{E\in \Th} \norm[E]{ \nabla u}^2}},
	\label{eq:errorH1}
\end{equation}
at increasing values of the local polynomial degree $k$ and at fixed mesh.

We expect an overall improvement with respect to the monomial approach, but not necessarily with respect to the orthonormal approach.
We recall that the aim of this work is to propose a computational cheap strategy to mitigate the ill-conditioning caused by the use of the standard scaled monomial basis defined on the original polytopal elements.
In this regard, we recall that the additional cost of the orthonormal approach relies on the cost of the application of the Modified Gram Schmidt algorithm with reorthogonalization on each element (and on each face, if $d=3$), which depends on $n_k^d = \frac{(k+1)\dots(k+d)}{d!}$.
On the contrary, the computational complexity of our method only depends on the geometric dimension of the problem $d$ and, if $d=3$, on the number of faces for the \ref{itm:F} and \ref{itm:BF} inertial approaches, but not on the local polynomial degree $k$.

\subsection{Test 1: Highly-distorted mesh, small domain and hanging nodes in 2D}
\label{sec:test1}
\begin{figure}[]
	\centering
	\subfigure[\label{fig:distortedhexa70x80}]
	{\includegraphics[width=.24\textwidth, height = .17\textheight]{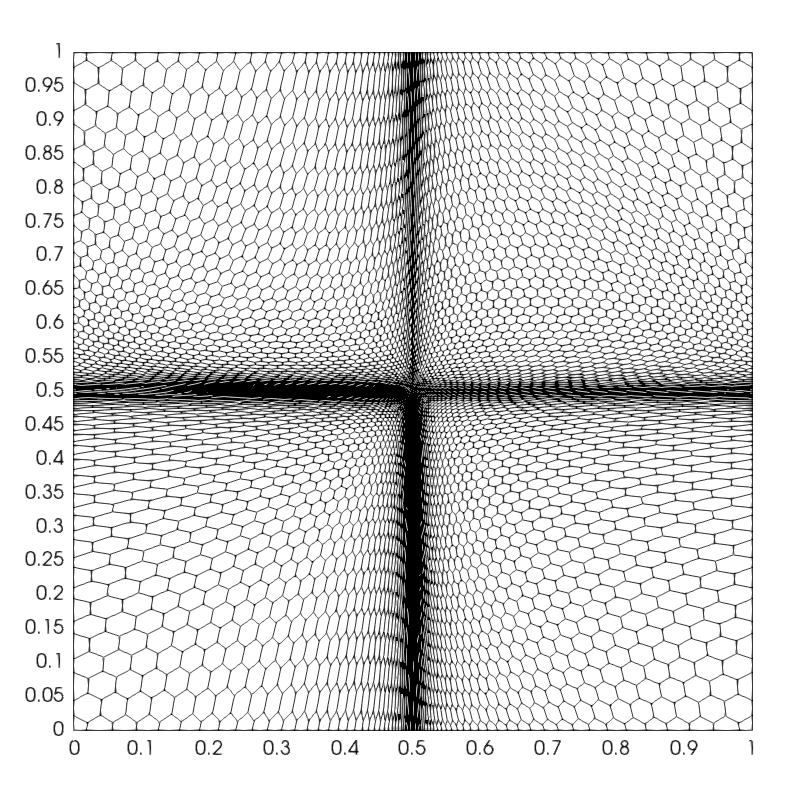}}
	\subfigure[\label{fig:tri_12}]
	{\includegraphics[width=.24\textwidth, height = .17\textheight]{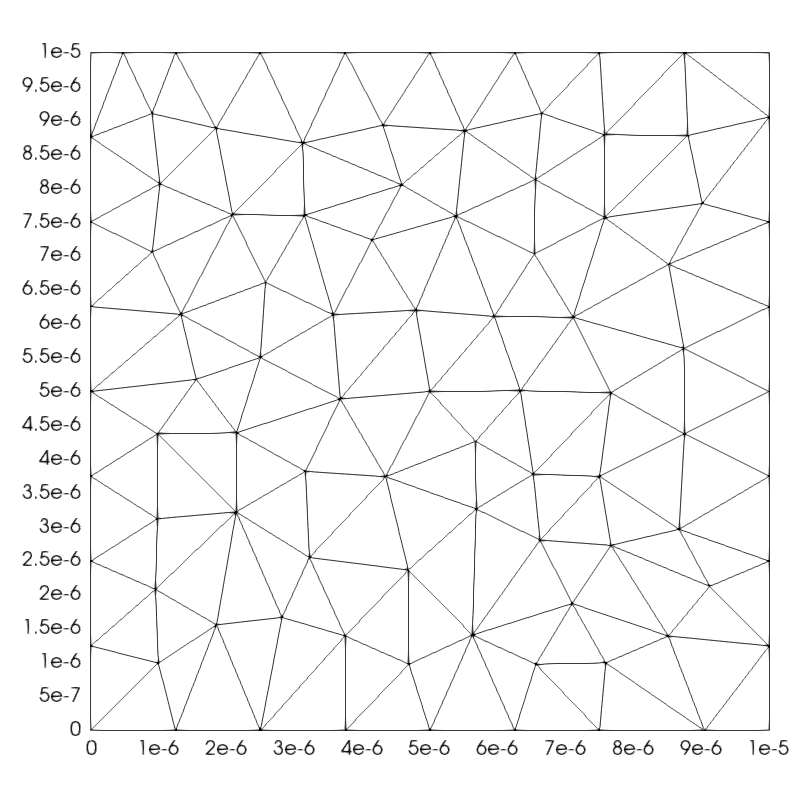}}
	\subfigure[\label{fig:polymesh2D}]
	{\includegraphics[width=.24\textwidth, height = .17\textheight]{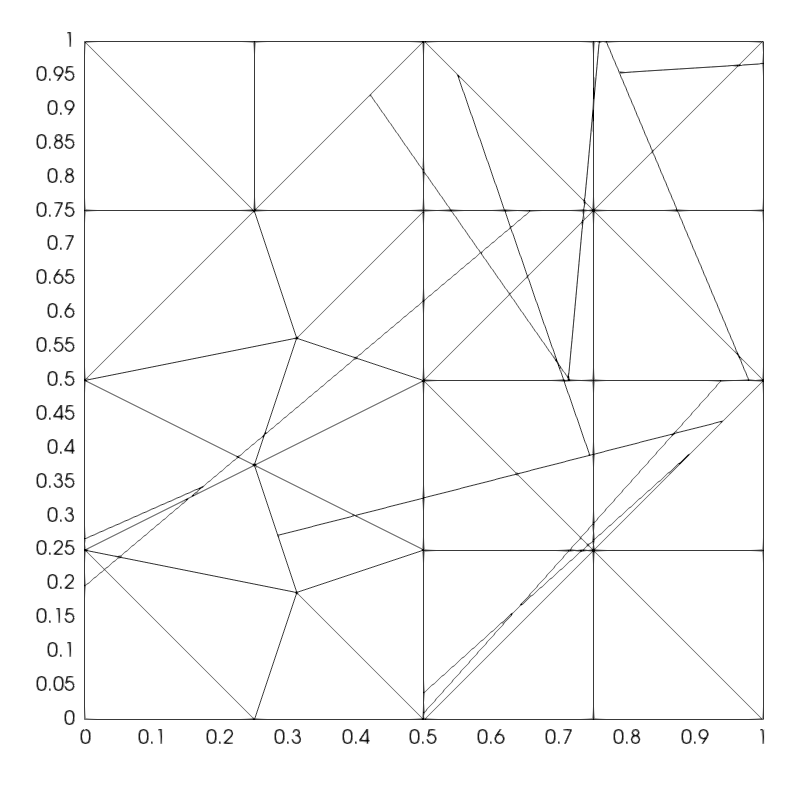}}
	\subfigure[\label{fig:Run_1e_2}]
	{\includegraphics[width=.24\textwidth, height = .17\textheight]{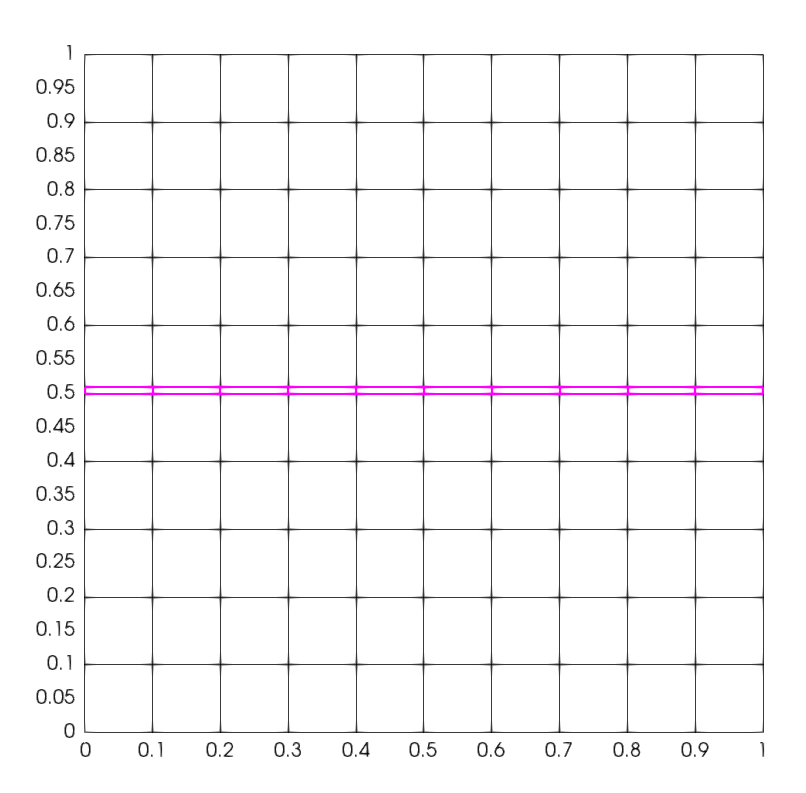}}
	\caption{Meshes used for tests 1 and 2. Highly-distorted hexagonal mesh (HDHM) on \ref{fig:distortedhexa70x80}. Triangular mesh on a small domain (RTRM) on \ref{fig:tri_12}. Polygonal mesh (GPGM) on \ref{fig:polymesh2D}. Collapsing polygons (CSM$_1$) on \ref{fig:Run_1e_2}.}
	\label{fig:mesh_2D}
\end{figure}
\begin{table}[]
	\centering
	\resizebox{\textwidth}{!}{
\begin{tabular}{lccccccccccccccc}
\multirow{2}{*}{}                 & \multirow{2}{*}{\textbf{$\#\mathcal{T}_h$}} & \textbf{$N_{v}^E$}   & \multirow{2}{*}{\textbf{}} & \multicolumn{3}{c|}{\textbf{Area}}                                                          & \multicolumn{3}{c|}{\textbf{Diameter}}                                                      & \multicolumn{3}{c|}{\textbf{Anisotropic  Ratio}}                                            & \multicolumn{3}{c|}{\textbf{Edge  Ratio}}                                                   \\
                                  &                                             & \textbf{avg}         &                            & \textbf{min}               & \textbf{avg}               & \multicolumn{1}{c|}{\textbf{max}} & \textbf{min}               & \textbf{avg}               & \multicolumn{1}{c|}{\textbf{max}} & \textbf{min}               & \textbf{avg}               & \multicolumn{1}{c|}{\textbf{max}} & \textbf{min}               & \textbf{avg}               & \multicolumn{1}{c|}{\textbf{max}} \\ \hline
\multirow{2}{*}{\textbf{HDHM}}    & \multirow{2}{*}{5711}                       & \multirow{2}{*}{5.9} & $E$                        & 2.6e-6                     & 1.8e-4                     & \multicolumn{1}{c|}{1.3e-3}       & 2.1e-3                     & 2.4e-2                     & \multicolumn{1}{c|}{6.2e-2}       & 1.0e+0                     & 5.9e+1                     & \multicolumn{1}{c|}{6.2e+2}       & 1.0e+0                     & 2.5e+0                     & \multicolumn{1}{c|}{2.0e+1}       \\
                                  &                                             &                      & $\hat{E}$                  & 5.0e-1                     & 6.4e-1                     & \multicolumn{1}{c|}{6.5e-1}       & 1.0e+0                     & 1.0e+0                     & \multicolumn{1}{c|}{1.0e+0}       & 1.0e+0                     & 1.0e+0                     & \multicolumn{1}{c|}{1.0e+0}       & 1.0e+0                     & 1.1e+0                     & \multicolumn{1}{c|}{2.0e+0}       \\ \hline
\multirow{2}{*}{\textbf{RTRM}}    & \multirow{2}{*}{153}                        & \multirow{2}{*}{3}   & $E$                        & 3.0e-13                    & 6.5e-13                    & \multicolumn{1}{c|}{1.0e-12}      & 1.0e-6                     & 1.4e-6                     & \multicolumn{1}{c|}{1.9e-6}       & 1.0e+0                     & 2.3e+0                     & \multicolumn{1}{c|}{9.9e+0}       & 1.0e+0                     & 1.3e+0                     & \multicolumn{1}{c|}{2.8e+0}       \\
                                  &                                             &                      & $\hat{E}$                  & 4.3e-1                     & 4.3e-1                     & \multicolumn{1}{c|}{4.3e-1}       & 1.0e+0                     & 1.0e+0                     & \multicolumn{1}{c|}{1.0e+0}       & 1.0e+0                     & 1.0e+0                     & \multicolumn{1}{c|}{1.0e+0}       & 1.0e+0                     & 1.0e+0                     & \multicolumn{1}{c|}{1.0e+0}       \\ \hline
\multirow{2}{*}{\textbf{GPGM}}    & \multirow{2}{*}{81}                         & \multirow{2}{*}{4.3} & $E$                        & 3.7e-6                     & 1.2e-2                     & \multicolumn{1}{c|}{3.1e-2}       & 3.7e-3                     & 2.0e-1                     & \multicolumn{1}{c|}{3.5e-1}       & 1.3e+0                     & 2.7e+1                     & \multicolumn{1}{c|}{4.7e+2}       & 1.2e+0                     & 5.1e+2                     & \multicolumn{1}{c|}{4.0e+4}       \\
                                  &                                             &                      & $\hat{E}$                  & \multicolumn{1}{l}{4.3e-1} & \multicolumn{1}{l}{4.5e-1} & \multicolumn{1}{l|}{5.3e-1}       & \multicolumn{1}{l}{1.0e+0} & \multicolumn{1}{l}{1.0e+0} & \multicolumn{1}{l|}{1.0e+0}       & \multicolumn{1}{l}{1.0e+0} & \multicolumn{1}{l}{1.0e+0} & \multicolumn{1}{l|}{1.0e+0}       & \multicolumn{1}{l}{1.0e+0} & \multicolumn{1}{l}{1.3e+2} & \multicolumn{1}{l|}{9.0e+3}       \\ \hline
                                  & \multicolumn{1}{l}{}                        & \multicolumn{1}{l}{} & \multicolumn{1}{l}{}       & \multicolumn{1}{l}{}       & \multicolumn{1}{l}{}       & \multicolumn{1}{l}{}              & \multicolumn{1}{l}{}       & \multicolumn{1}{l}{}       & \multicolumn{1}{l}{}              & \multicolumn{1}{l}{}       & \multicolumn{1}{l}{}       & \multicolumn{1}{l}{}              & \multicolumn{1}{l}{}       & \multicolumn{1}{l}{}       & \multicolumn{1}{l}{}              \\ \hline
\multirow{2}{*}{\textbf{CSM$_1$}} & \multirow{2}{*}{110}                        & \multirow{2}{*}{4}   & $E$                        & 1.0e-3                     & 9.1e-3                     & \multicolumn{1}{c|}{1.0e-2}       & 1.0e-1                     & 1.4e-1                     & \multicolumn{1}{c|}{1.4e-1}       & 1.0e+0                     & 1.0e+1                     & \multicolumn{1}{c|}{1.0e+2}       & 1.0e+0                     & 1.8e+0                     & \multicolumn{1}{c|}{1.0e+1}       \\
                                  &                                             &                      & $\hat{E}$                  & 5.0e-1                     & 5.0e-1                     & \multicolumn{1}{c|}{5.0e-1}       & 1.0e+0                     & 1.0e+0                     & \multicolumn{1}{c|}{1.0e+0}       & 1.0e+0                     & 1.0e+0                     & \multicolumn{1}{c|}{1.0e+0}       & 1.0e+0                     & 1.0e+0                     & \multicolumn{1}{c|}{1.0e+0}       \\ \hline
\multirow{2}{*}{\textbf{CSM$_2$}} & \multirow{2}{*}{110}                        & \multirow{2}{*}{4}   & $E$                        & 1.0e-4                     & 9.1e-3                     & \multicolumn{1}{c|}{1.0e-2}       & 1.0e-1                     & 1.4e-1                     & \multicolumn{1}{c|}{1.4e-1}       & 1.0e+0                     & 9.1e+2                     & \multicolumn{1}{c|}{1.0e+4}       & 1.0e+0                     & 1.0e+1                     & \multicolumn{1}{c|}{1.0e+2}       \\
                                  &                                             &                      & $\hat{E}$                  & 5.0e-1                     & 5.0e-1                     & \multicolumn{1}{c|}{5.0e-1}       & 1.0e+0                     & 1.0e+0                     & \multicolumn{1}{c|}{1.0e+0}       & 1.0e+0                     & 1.0e+0                     & \multicolumn{1}{c|}{1.0e+0}       & 1.0e+0                     & 1.0e+0                     & \multicolumn{1}{c|}{1.0e+0}       \\ \hline
\multirow{2}{*}{\textbf{CSM$_3$}} & \multirow{2}{*}{110}                        & \multirow{2}{*}{4}   & $E$                        & 1.0e-5                     & 9.1e-3                     & \multicolumn{1}{c|}{1.0e-2}       & 1.0e-1                     & 1.4e-1                     & \multicolumn{1}{c|}{1.4e-1}       & 1.0e+0                     & 9.1e+4                     & \multicolumn{1}{c|}{1.0e+6}       & 1.0e+0                     & 9.2e+1                     & \multicolumn{1}{c|}{1.0e+3}       \\
                                  &                                             &                      & $\hat{E}$                  & 5.0e-1                     & 5.0e-1                     & \multicolumn{1}{c|}{5.0e-1}       & 1.0e+0                     & 1.0e+0                     & \multicolumn{1}{c|}{1.0e+0}       & 1.0e+0                     & 1.0e+0                     & \multicolumn{1}{c|}{1.0e+0}       & 1.0e+0                     & 1.0e+0                     & \multicolumn{1}{c|}{1.0e+0}       \\ \hline
\end{tabular}}
	\caption{Properties of meshes used in the Test 1 and in the Test 2.}
	\label{tab:properties_mesh2D}
\end{table}

\begin{figure}[]
	\centering
	\subfigure[\label{fig:condPi0Derx_distortedHexa70x80_EpsilonPoisson_1.000000}]
	{\includegraphics[width=.32\textwidth, height = .17\textheight]{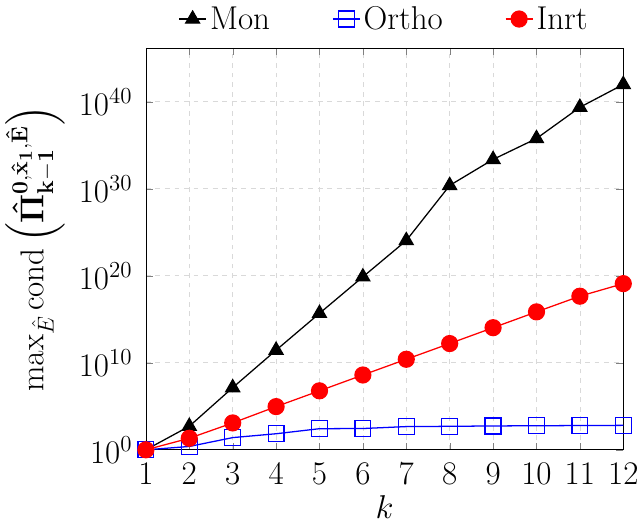}}
	\subfigure[\label{fig:condPi0Derx_tri_12_EpsilonPoisson_0.000010}]
	{\includegraphics[width=.32\textwidth, height = .17\textheight]{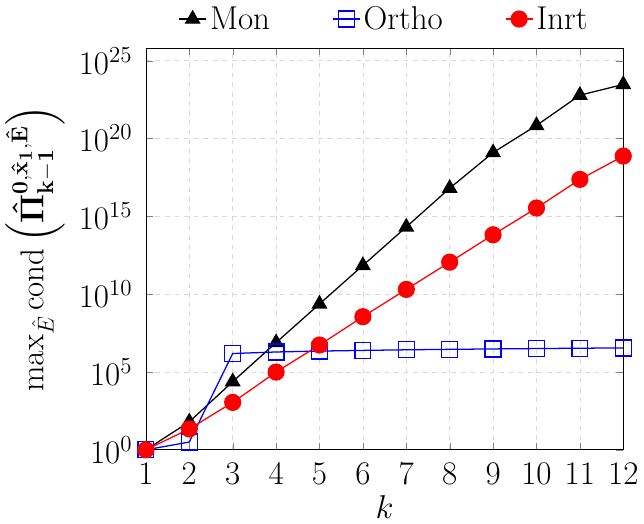}} 	\subfigure[\label{fig:condPi0Derx_polymesh2D_EpsilonPoisson_1.000000}]
	{\includegraphics[width=.32\textwidth, height = .17\textheight]{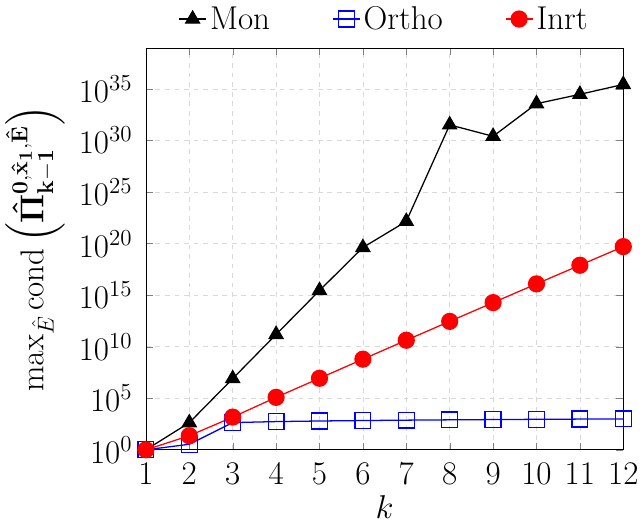}}

	\subfigure[\label{fig:condPiNabla_distortedHexa70x80_EpsilonPoisson_1.000000}]{\includegraphics[width=.32\textwidth, height = .17\textheight]{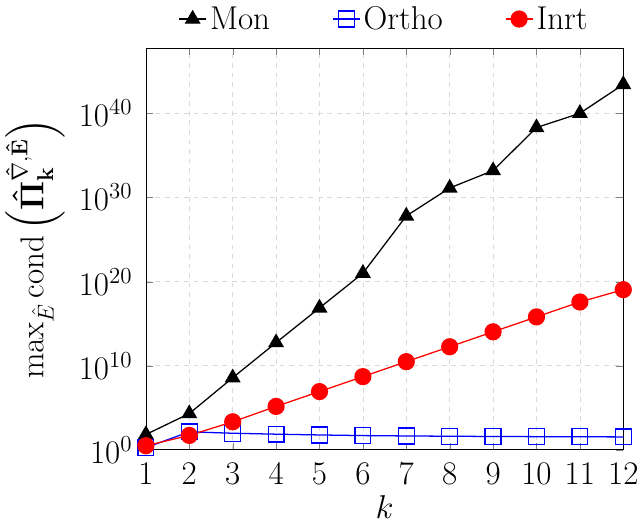}}
	\subfigure[\label{fig:condPiNabla_tri_12_EpsilonPoisson_0.000010}]
	{\includegraphics[width=.32\textwidth, height = .17\textheight]{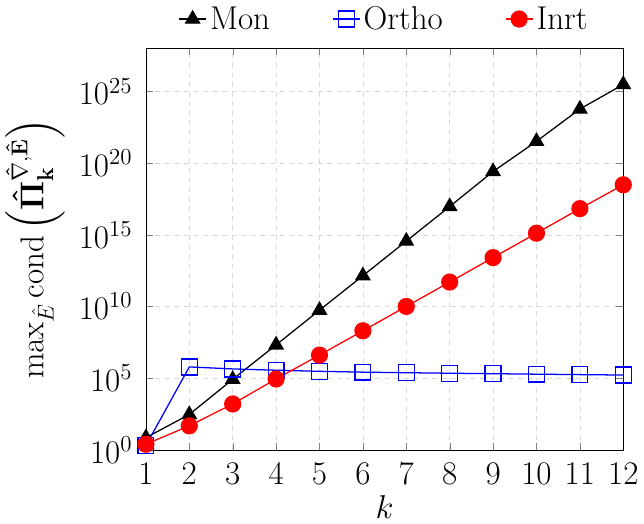}}	\subfigure[\label{fig:condPiNabla_polymesh2D_EpsilonPoisson_1.000000}]
	{\includegraphics[width=.32\textwidth, height = .17\textheight]{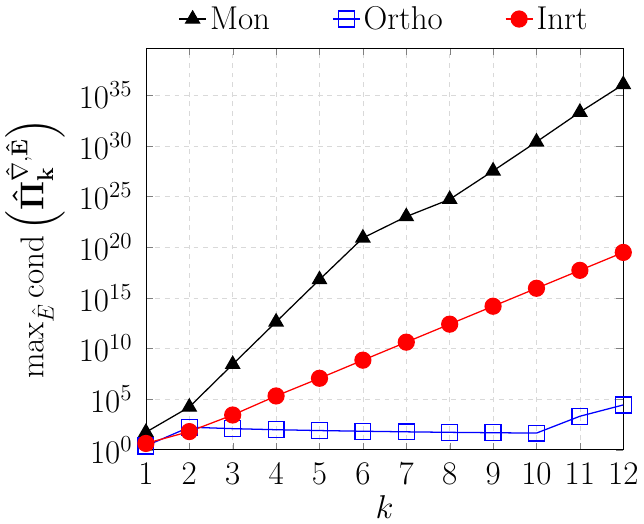}}

	\subfigure[\label{fig:condPi0km1_distortedHexa70x80_EpsilonPoisson_1.000000}]{\includegraphics[width=.32\textwidth, height = .17\textheight]{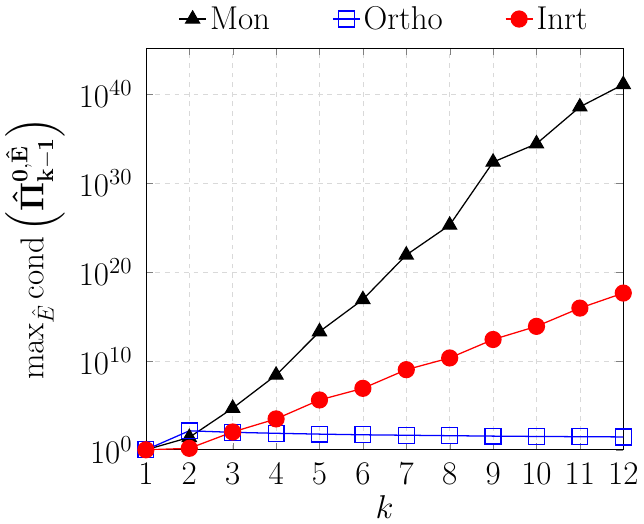}}
	\subfigure[\label{fig:condPi0km1_tri_12_EpsilonPoisson_0.000010}]
	{\includegraphics[width=.32\textwidth, height = .17\textheight]{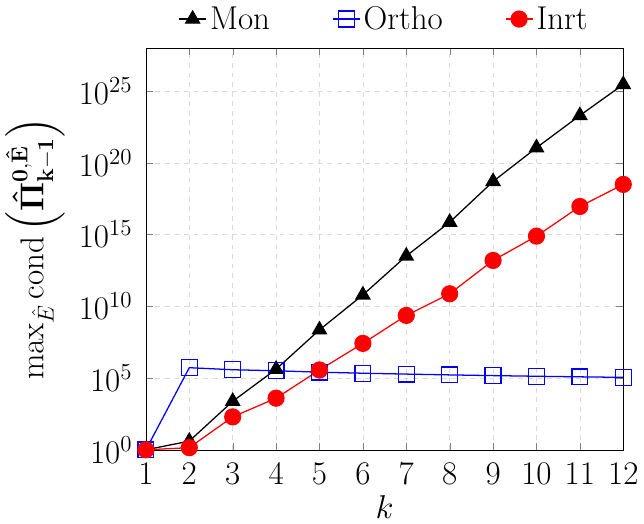}}	\subfigure[\label{fig:condPi0km1_polymesh2D_EpsilonPoisson_1.000000}]
	{\includegraphics[width=.32\textwidth, height = .17\textheight]{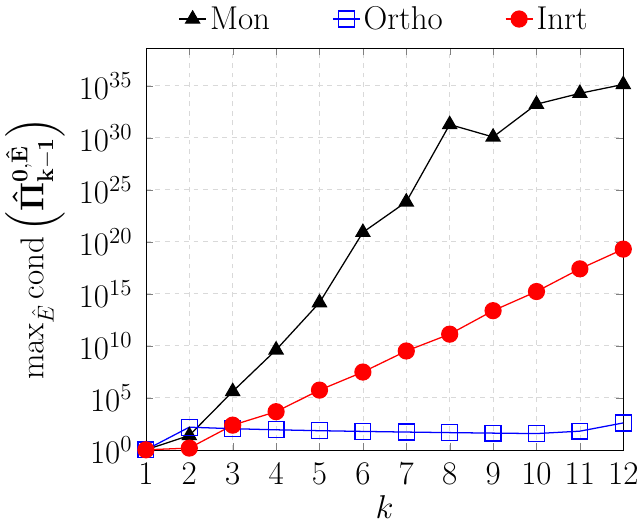}}

	\caption{Test 1: Behaviours of the worst condition numbers of local projection matrices among elements at varying $k$. Left: HDHM. Center: RTRM. Right: GPGM.}
	\label{fig:condlocalprojection_test1}
\end{figure}
\begin{figure}[]
	\centering
	\subfigure[\label{fig:condStiff_distortedHexa70x80_EpsilonPoisson_1.000000}]
	{\includegraphics[width=.32\textwidth, height = .17\textheight]{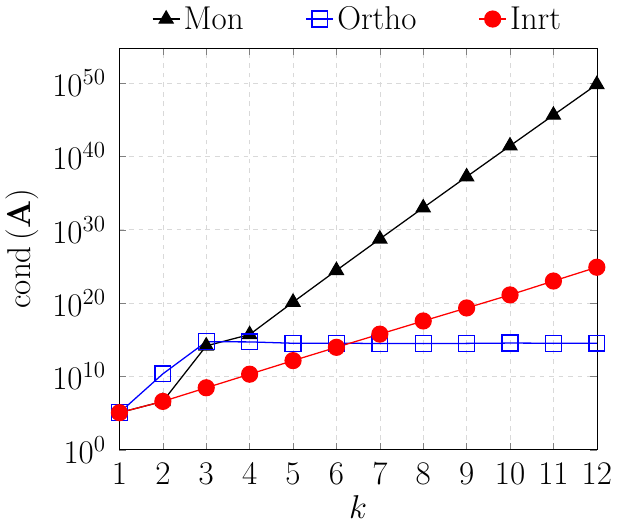}}
	\subfigure[\label{fig:condStiff_tri_12_EpsilonPoisson_0.000010}]
	{\includegraphics[width=.32\textwidth, height = .17\textheight]{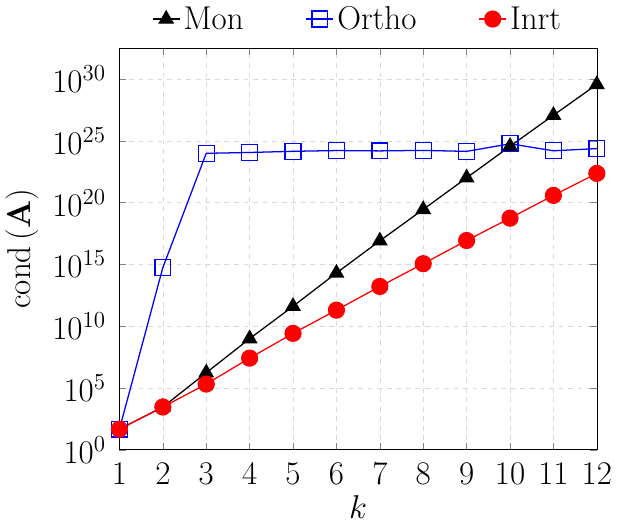}}
	\subfigure[\label{fig:condStiff_polymesh2D_EpsilonPoisson_1.000000}]
	{\includegraphics[width=.32\textwidth, height = .17\textheight]{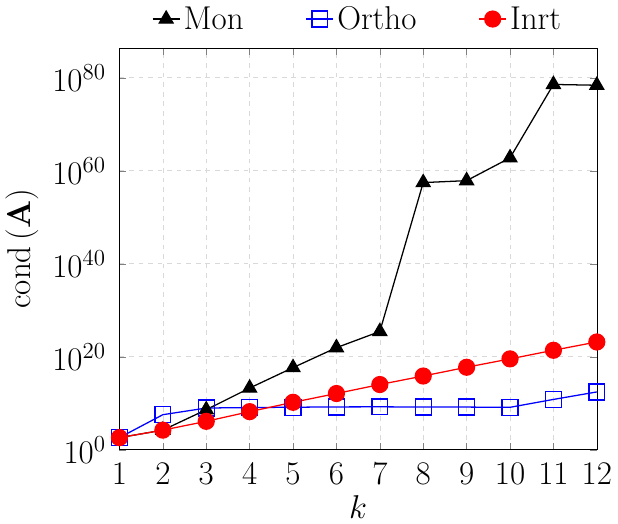}}
	\caption{Test 1: Behaviour of the condition number of the system matrix $\mathbf{A}$ at varying $k$. Left: HDHM. Center: RTRM. Right: GPGM.}
	\label{fig:condstiff_test1}
\end{figure}
\begin{figure}[]
	\centering
	\subfigure[\label{fig:ErrorL2_distortedHexa70x80_EpsilonPoisson_1.000000}]
	{\includegraphics[width=.32\textwidth, height = .17\textheight]{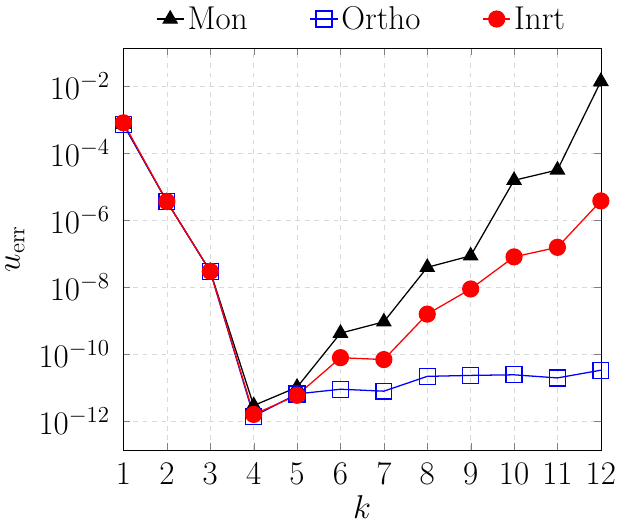}}
	\subfigure[\label{fig:ErrorL2_tri_12_EpsilonPoisson_0.000010}]
	{\includegraphics[width=.32\textwidth, height = .17\textheight]{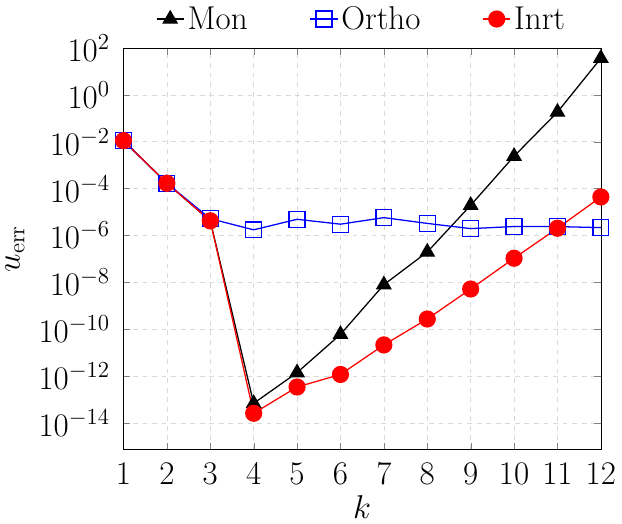}}
	\subfigure[\label{fig:ErrorL2_polymesh2D_EpsilonPoisson_1.000000}]
	{\includegraphics[width=.32\textwidth, height = .17\textheight]{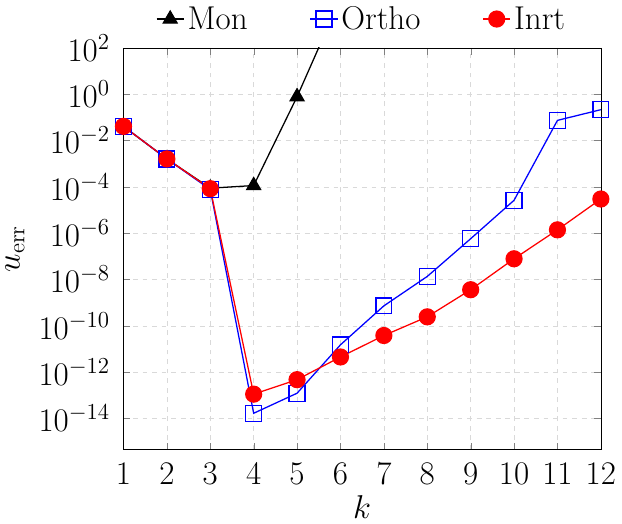}}
	\subfigure[\label{fig:ErrorH1_distortedHexa70x80_EpsilonPoisson_1.000000}]
	{\includegraphics[width=.32\textwidth, height = .17\textheight]{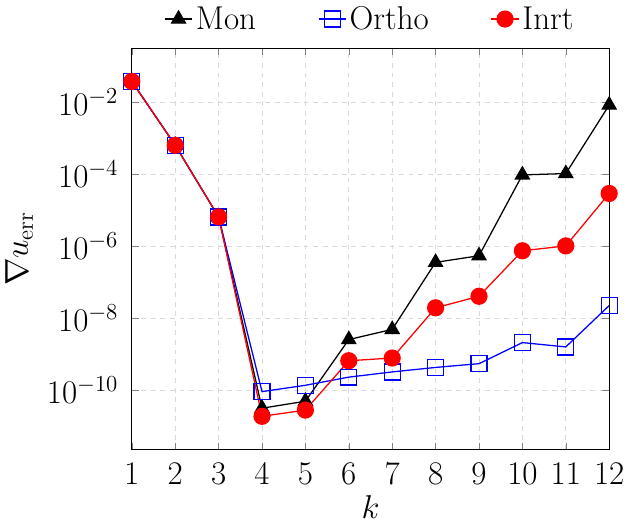}}
	\subfigure[\label{fig:ErrorH1_tri_12_EpsilonPoisson_0.000010}]
	{\includegraphics[width=.32\textwidth, height = .17\textheight]{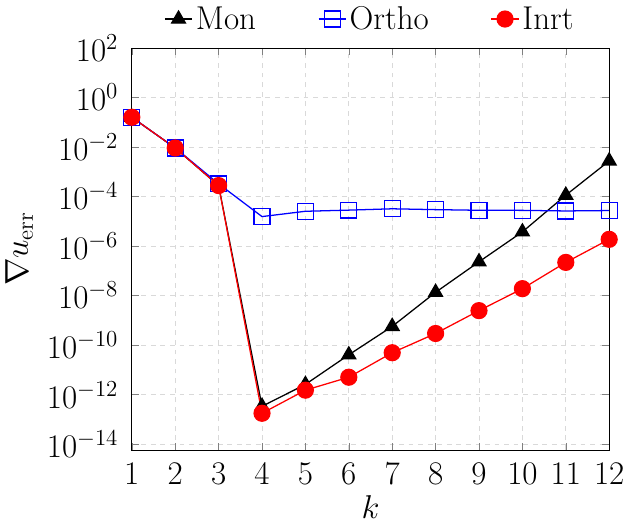}}
	\subfigure[\label{fig:ErrorH1_polymesh2D_EpsilonPoisson_1.000000}]
	{\includegraphics[width=.32\textwidth, height = .17\textheight]{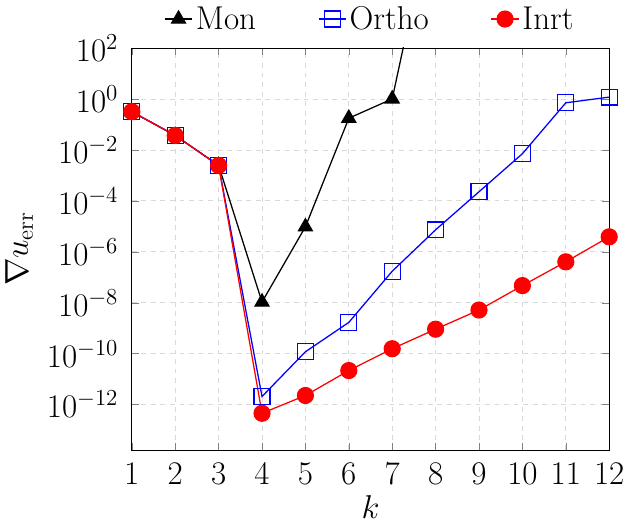}}
	\caption{Test 1: Behaviours of errors \eqref{eq:errorL2} and \eqref{eq:errorH1} at varying $k$. Left: HDHM. Center: RTRM. Right: GPGM.}
	\label{fig:error_test1}
\end{figure}

Given $\epsilon \in \R^{+}$ and $\Omega =\left(0,\epsilon\right)\times \left(0,\epsilon\right)$, we consider the problem \eqref{eq:modelproblem} with constant coefficients $\D = \mathbf{I}$, $\gamma = 0$ and $\bb= \begin{bmatrix}
0\\0
\end{bmatrix}$, where $f$ and the non-homogeneous Dirichlet boundary condition are set in such a way the exact solution is:

\begin{equation}
	u(x_1,x_2) = 1.1 + \frac{16}{\epsilon^4} x_1x_2(\epsilon-x_1)(\epsilon-x_2).
	\label{eq:solution_test1}
\end{equation}
In this first test, we consider
\begin{itemize}
	\item a highly-distorted hexagonal mesh (HDHM in short) on a square domain with edge length $\epsilon = 1$ (Figure~\ref{fig:distortedhexa70x80});
	\item a regular triangular mesh (RTRM in short) on a small square domain with edge length $\epsilon = 1.0e\text{-}5$ (Figure~\ref{fig:tri_12});
	\item a generic polygonal mesh (GPGM in short) on a square domain with edge length $\epsilon = 1$, which is characterised by polygons with very different shapes and areas and by the presence of hanging nodes (Figure~\ref{fig:polymesh2D}).
\end{itemize}

In Table~\ref{tab:properties_mesh2D}, we show the main features of both the original and the mapped polygons related to each aforementioned mesh, namely the area, the diameter, the anisotropic ratio and the \emph{edge ratio}, i.e. the ratio between the highest and the smallest lengths of edges of $E$. 
By looking at the geometric properties of $\srescale{E}$, we note that the map $F_E$ generates well-shaped polygons accordingly to our definition. Furthermore, it tends to uniform the geometric properties of elements belonging to the same category (such as triangles, parallelograms, hexagons,\dots) in absence of almost-hanging nodes, i.e. nodes between two consecutive edges forming an angle of about 180 degrees.
Indeed, in presence of almost-hanging nodes, our map does not eliminate any problems related to small edges that participate to form almost-hanging nodes, as highlighted by looking at the edge ratio property of the mesh GPGM. 
 
In figures ~\ref{fig:condlocalprojection_test1}, \ref{fig:condstiff_test1} and ~\ref{fig:error_test1}, we show the trends of the condition numbers of the local projection matrices, of the global system matrix and of the errors \eqref{eq:errorL2}-\eqref{eq:errorH1} at varying the polynomial degree $k$, respectively. We omit the graphs reporting the behaviours of the condition numbers of $\mmproj{0,\srescale{x}_2,\srescale{E}}{k-1}$ at varying $k$, since the trend is very similar to the one of the condition number of $\mmproj{0,\srescale{x}_1,\srescale{E}}{k-1}$ for each method.
In all the figures, we note that our proposed procedure outperforms the monomial one as expected.

In Figure~\ref{fig:condlocalprojection_test1}, we observe good local results of the inertial approach, despite the best performances are obtained by the orthonormal approach for high values of $k$. 
To this regard, we remark that our aim is to cheaply reduce the ill-conditioning of local and global matrices in presence of badly-shaped elements with respect to the monomial approach.

We stress that, in the case of mesh RTRM, the first re-scaling in \eqref{eq:mapping} is required for the Inrt approach. Indeed, the mass matrix related to the original RTRM elements is a singular matrix in finite precision due to its eigenvalues which are in the order of magnitude of the round-off error. We point out that the small triangles of RTRM represent a challenging geometry for the orthonormal approach. The inertial approach, instead, is robust in presence of the small polygons of mesh RTRM in terms of the condition number of the global matrix (see Figure ~\ref{fig:condStiff_tri_12_EpsilonPoisson_0.000010}). 

Moreover, in presence of very badly-shaped polygons of mesh GPGM, the global performances of the Ortho and Inrt methods are comparable (see Figure \ref{fig:condStiff_polymesh2D_EpsilonPoisson_1.000000}).

Since the proposed solution \eqref{eq:solution_test1} is a polynomial of degree 4, we expect that the errors shown in Figure~\ref{fig:error_test1} tend to zero when $k<4$ and then vanish when $k \geq 4$. However, after an initial decrease, in most cases, the errors start to raise due to the ill-conditioning. In this regard, we note that the monomial errors blow up in the case of mesh GPGM, 
while the growth of errors related to Inrt is much more controlled. 
Finally, we remark that the robustness of our method in the cases of RTRM and GPGM meshes reflects in the smallest errors for the higher values of $k$, as highlighted in figures \ref{fig:ErrorL2_tri_12_EpsilonPoisson_0.000010},  \ref{fig:ErrorL2_polymesh2D_EpsilonPoisson_1.000000},\ref{fig:ErrorH1_tri_12_EpsilonPoisson_0.000010} and \ref{fig:ErrorH1_polymesh2D_EpsilonPoisson_1.000000}.

\subsection{Test 2: Collapsing polygons}
\begin{figure}[]
	\centering
	\subfigure[\label{fig:condStifcondPi0Derx_Run_Ellipticf_Run_Elliptic}]
	{\includegraphics[width=.32\textwidth, height = .17\textheight]{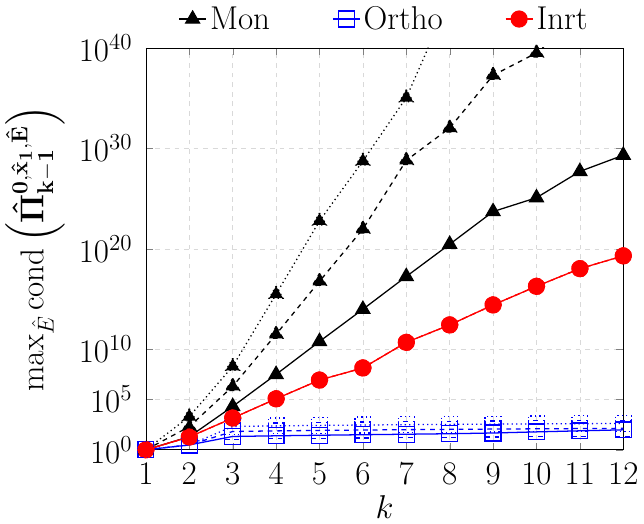}}
	\subfigure[\label{fig:condPiNabla_Run_Elliptic}]
	{\includegraphics[width=.32\textwidth, height = .17\textheight]{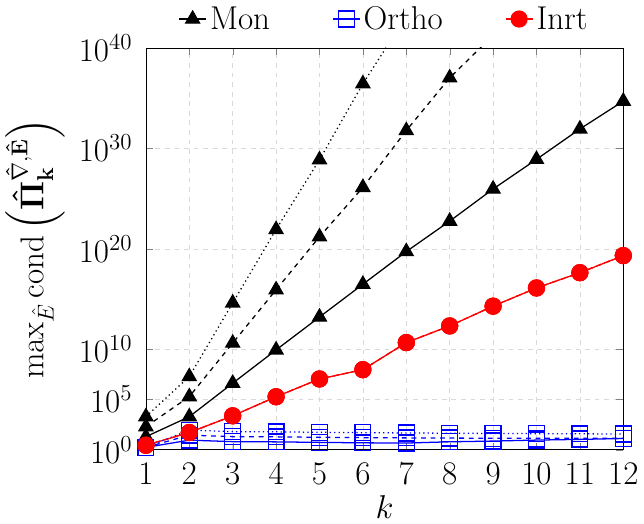}}
	\subfigure[\label{fig:condPi0km1_Run_Elliptic}]
	{\includegraphics[width=.32\textwidth, height = .17\textheight]{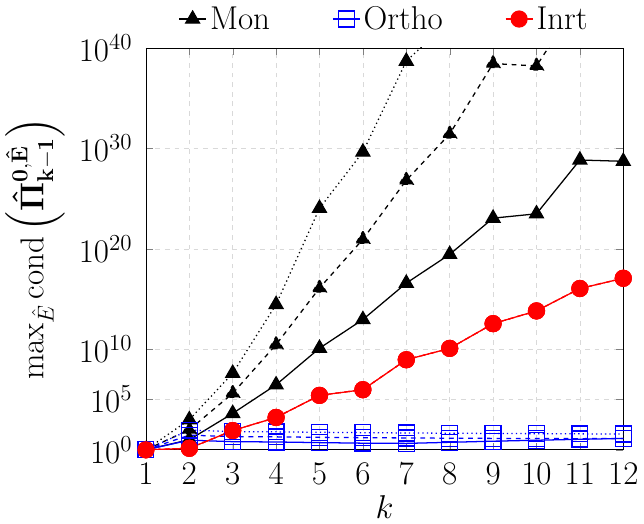}}
	\caption{Test 2: Behaviours of the worst condition numbers of local projection matrices among elements at varying $k$. Solid lines: CSM$_1$. Dashed lines: CSM$_2$. Dotted lines: CSM$_3$.}
	\label{fig:condlocalprojection_test2}
\end{figure}
\begin{figure}[]
	\centering
	\subfigure[\label{fig:condStiff_Run_Elliptic}]
	{\includegraphics[width=.32\textwidth, height = .17\textheight]{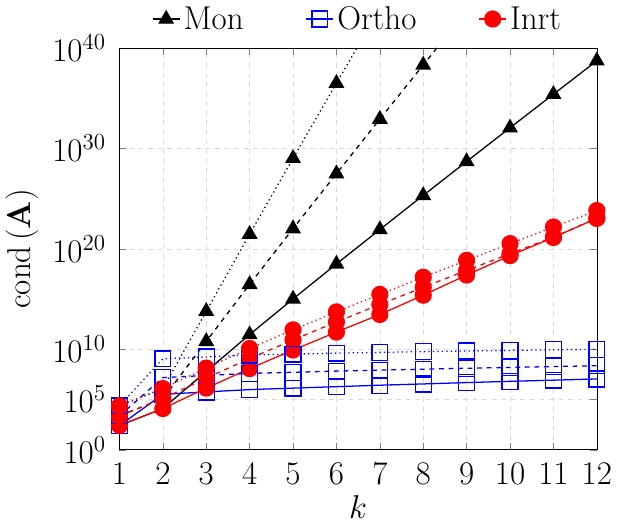}}
	\subfigure[\label{fig:ErrorL2_Run_Elliptic}]
	{\includegraphics[width=.32\textwidth, height = .17\textheight]{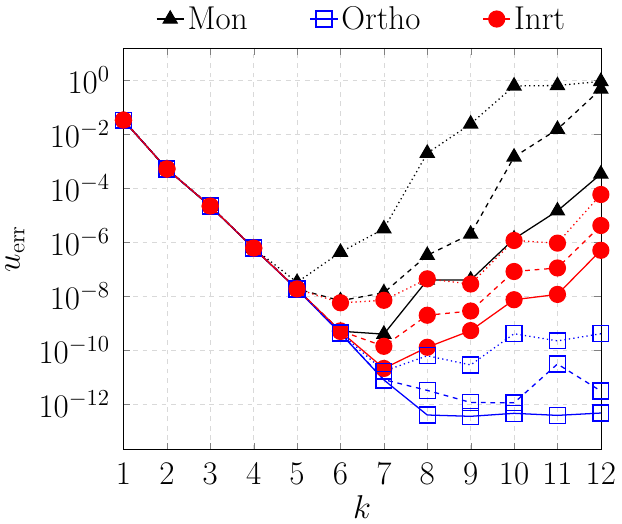}}
	\subfigure[\label{fig:ErrorH1_Run_Elliptic}]
	{\includegraphics[width=.32\textwidth, height = .17\textheight]{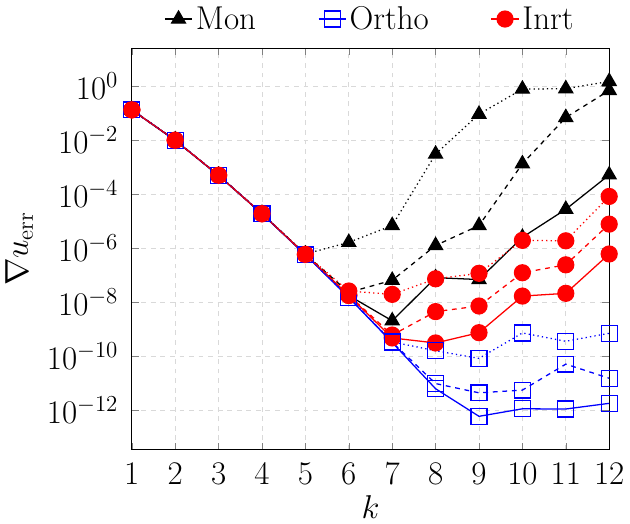}}
	\caption{Test 2: Behaviours of the condition numbers of the global system matrix and of errors \eqref{eq:errorL2}, \eqref{eq:errorH1} at varying $k$. Solid lines: CSM$_1$. Dashed lines: CSM$_2$. Dotted lines: CSM$_3$.}
	\label{fig:error_condstiff_test2}
\end{figure}
In this experiment, we test the performances of our procedures considering the complete diffusion-advection-reaction problem \eqref{eq:modelproblem} with variable coefficients and a non-polynomial solution in the case of collapsing polygons.
Thus, on domain $\Omega =\left(0,1\right)\times \left(0,1\right)$, we consider

\begin{equation*}
	\D(x_1,x_2) = \begin{bmatrix}
		1 + x_2^2 & -x_1x_2 \\
		-x_1x_2   & 1+x_1^2
	\end{bmatrix},\quad \bb(x_1,x_2) = \begin{bmatrix}
		x_1 \\ -x_2
	\end{bmatrix}, \quad \gamma(x_1,x_2) = x_1x_2,
\end{equation*}
and we set $f$ in such a way the exact solution is

\begin{equation*}
	u(x_1,x_2) = \sin(\pi x_1) \sin(\pi x_2).
\end{equation*}

We generate on the domain a sequence of three rectangular meshes composed by squared elements of area $10^{-1}$ with the exception of a central band made up of two groups of rectangles,
one of which (the purple band highlighted in Figure~\ref{fig:Run_1e_2}) is formed by rectangles of height $10^{-2}$, $10^{-3}$ and $10^{-4}$ for CSM$_1$, CSM$_2$ and CSM$_3$ mesh respectively. We refer to Figure~\ref{fig:Run_1e_2} for the plot of the first mesh CSM$_1$.

The geometric properties of the original and of the mapped elements of the three meshes are shown in Table~\ref{tab:properties_mesh2D} on the rows CSM$_i$, $i \in \{1,2,3\}$. 
We note that the performed map makes the mapped elements belonging to the central band equal to the other mapped ones. 
This result is well highlighted by looking at the condition numbers of local projection matrices in Figure~\ref{fig:condlocalprojection_test2}, which do not vary among the three meshes in the inertial approach. 
From the condition number of the global system matrix and from the error measurements of Figure~\ref{fig:error_condstiff_test2}, we can observe that the global performances remain dependent of the geometric properties of the original elements also in the inertial approach. As done for Test 1, we omit to report the behaviours of the condition numbers of $\mmproj{0,\srescale{x}_2,\srescale{E}}{k-1}$, since its trend is very similar to the one of the condition number of $\mmproj{0,\srescale{x}_1,\srescale{E}}{k-1}$ for each method. 


Finally, the condition and error measurements confirm that our technique shows very good improvements with respect to the monomial case. 
For high polynomial degrees, the orthonormal case outperforms the inertial one. However, we recall that our procedure leads to a huge reduction in the overall computational cost as it does not depend on the local polynomial degree $k$.

\subsection{Test 3: Badly-shaped polyhedrons and aligned faces in 3D}
\begin{figure}[]
	\centering
	\subfigure[\label{fig:mesh_Tetra200}]
	{\includegraphics[width=.32\textwidth, height = .17\textheight]{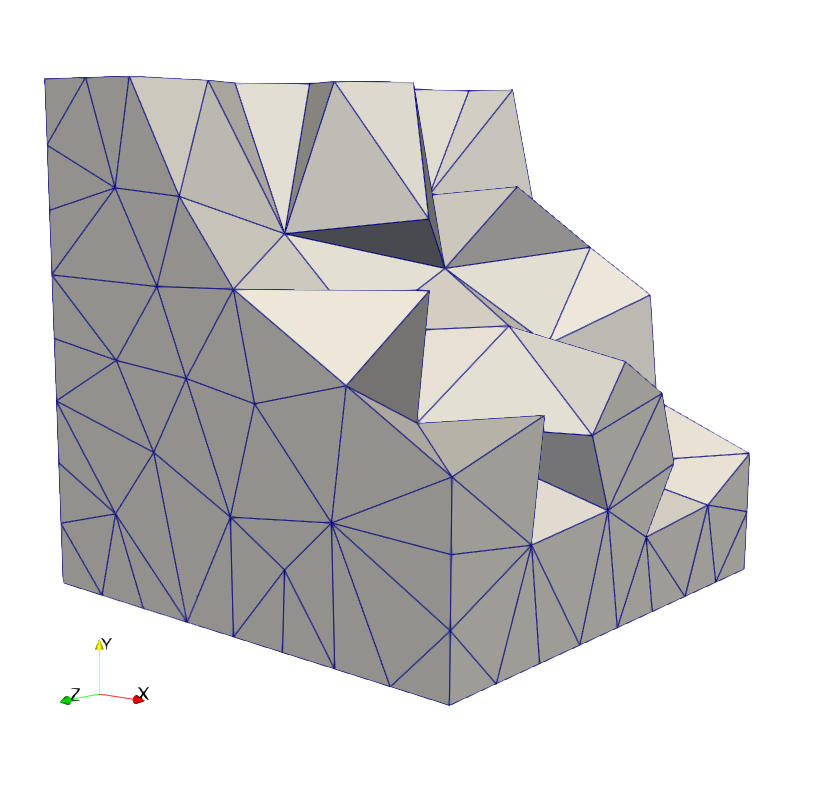}}
	\subfigure[\label{fig:mesh_conf}]
	{\includegraphics[width=.32\textwidth, height = .17\textheight]{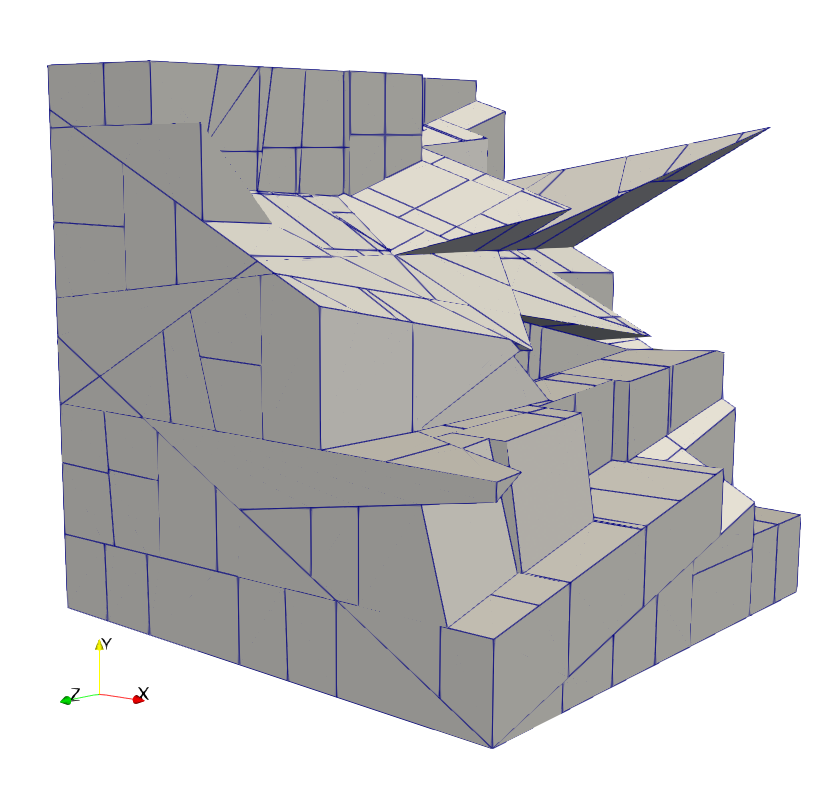}}
	\subfigure[\label{fig:mesh_Hexahedron_1e_01}]
	{\includegraphics[width=.32\textwidth, height = .17\textheight]{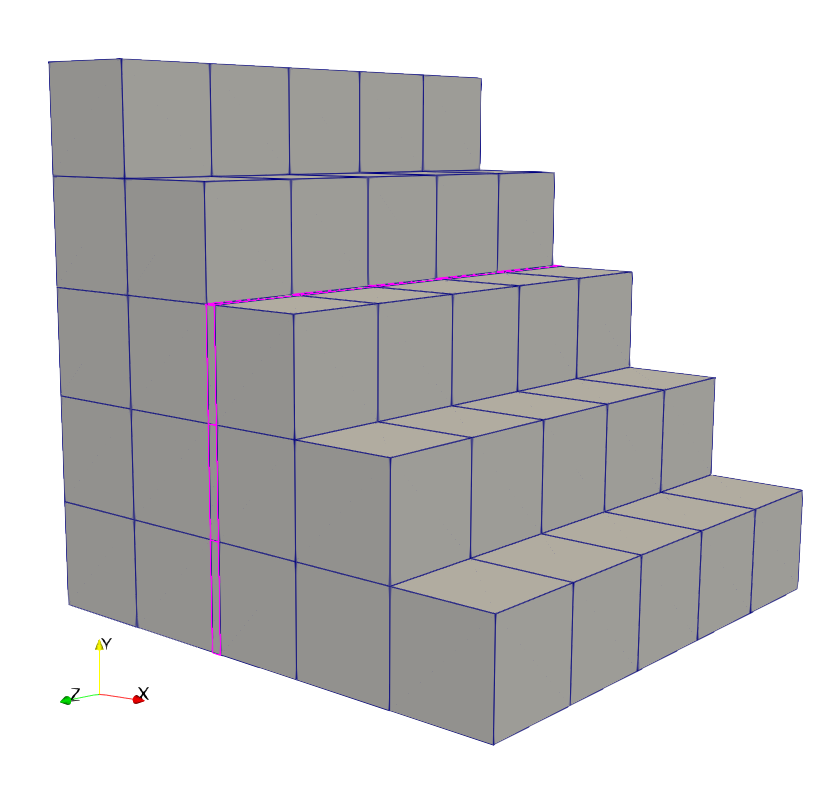}}
	\caption{Meshes used for Tests 3 and 4 clipped on the cubic domain $[0, 1]^3$. Left: Tetrahedral mesh (RTTM). Center: Polyhedral mesh (GPDM). Right: Collapsing polyhedrons (CCM$_1$).}
	\label{fig:mesh_3D}
\end{figure}
\begin{table}[]
	\centering
	\resizebox{\textwidth}{!}{
\begin{tabular}{lccccccccccccccc}
\multirow{2}{*}{}                 & \multirow{2}{*}{\textbf{$\#\mathcal{F}_h$}} & \textbf{$N_{v}^f$}   & \multirow{2}{*}{\textbf{}} & \multicolumn{3}{c|}{\textbf{Area}}                                              & \multicolumn{3}{c|}{\textbf{Diameter}}                                          & \multicolumn{3}{c|}{\textbf{Anisotropic  Ratio}}                                & \multicolumn{3}{c|}{\textbf{Edge  Ratio}}                                       \\
                                  &                                             & \textbf{avg}         &                            & \textbf{min}         & \textbf{avg}         & \multicolumn{1}{c|}{\textbf{max}} & \textbf{min}         & \textbf{avg}         & \multicolumn{1}{c|}{\textbf{max}} & \textbf{min}         & \textbf{avg}         & \multicolumn{1}{c|}{\textbf{max}} & \textbf{min}         & \textbf{avg}         & \multicolumn{1}{c|}{\textbf{max}} \\ \hline
\multirow{2}{*}{\textbf{RTTM}}    & \multirow{2}{*}{1330}                       & \multirow{2}{*}{3}   & $f$                        & 7.8e-3               & 3.0e-2               & \multicolumn{1}{c|}{6.8e-2}       & 1.8e-1               & 3.3e-1               & \multicolumn{1}{c|}{4.7e-1}       & 1.0e+0               & 3.4e+0               & \multicolumn{1}{c|}{1.2e+1}       & 1.0e+0               & 1.6e+0               & \multicolumn{1}{c|}{2.8e+0}       \\
                                  &                                             &                      & $\check{f}$                & 4.3e-1               & 4.3e-1               & \multicolumn{1}{c|}{4.3e-1}       & 1.0e+0               & 1.0e+0               & \multicolumn{1}{c|}{1.0e+0}       & 1.0e+0               & 1.0e+0               & \multicolumn{1}{c|}{1.0e+0}       & 1.0e+0               & 1.0e+0               & \multicolumn{1}{c|}{1.0e+0}       \\ \hline
\multirow{2}{*}{\textbf{GPDM}}    & \multirow{2}{*}{2017}                       & \multirow{2}{*}{5.1} & $f$                        & 3.6e-7               & 1.2e-2               & \multicolumn{1}{c|}{8.3e-2}       & 1.9e-3               & 1.7e-1               & \multicolumn{1}{c|}{7.9e-1}       & 1.0e+0               & 3.5e+3               & \multicolumn{1}{c|}{1.2e+6}       & 1.0e+0               & 4.4e+1               & \multicolumn{1}{c|}{4.4e+3}       \\
                                  &                                             &                      & $\check{f}$                & 4.3e-1               & 4.9e-1               & \multicolumn{1}{c|}{6.1e-1}       & 1.0e+0               & 1.0e+0               & \multicolumn{1}{c|}{1.0e+0}       & 1.0e+0               & 1.0e+0               & \multicolumn{1}{c|}{1.0e+0}       & 1.0e+0               & 2.6e+1               & \multicolumn{1}{c|}{2.7e+3}       \\ \hline
                                  & \multicolumn{1}{l}{}                        & \multicolumn{1}{l}{} & \multicolumn{1}{l}{}       & \multicolumn{1}{l}{} & \multicolumn{1}{l}{} & \multicolumn{1}{l}{}              & \multicolumn{1}{l}{} & \multicolumn{1}{l}{} & \multicolumn{1}{l}{}              & \multicolumn{1}{l}{} & \multicolumn{1}{l}{} & \multicolumn{1}{l}{}              & \multicolumn{1}{l}{} & \multicolumn{1}{l}{} & \multicolumn{1}{l}{}              \\ \hline
\multirow{2}{*}{\textbf{CCM$_1$}} & \multirow{2}{*}{535}                        & \multirow{2}{*}{4}   & $f$                        & 4.0e-3               & 3.6e-2               & \multicolumn{1}{c|}{4.0e-2}       & 2.0e-1               & 2.7e-1               & \multicolumn{1}{c|}{2.8e-1}       & 1.0e+0               & 1.2e+1               & \multicolumn{1}{c|}{1.0e+2}       & 1.0e+0               & 2.0e+0               & \multicolumn{1}{c|}{1.0e+1}       \\
                                  &                                             &                      & $\check{f}$                & 5.0e-1               & 5.0e-1               & \multicolumn{1}{c|}{5.0e-1}       & 1.0e+0               & 1.0e+0               & \multicolumn{1}{c|}{1.0e+0}       & 1.0e+0               & 1.0e+0               & \multicolumn{1}{c|}{1.0e+0}       & 1.0e+0               & 1.0e+0               & \multicolumn{1}{c|}{1.0e+0}       \\ \hline
\multirow{2}{*}{\textbf{CCM$_2$}} & \multirow{2}{*}{535}                        & \multirow{2}{*}{4}   & $f$                        & 4.0e-4               & 3.6e-2               & \multicolumn{1}{c|}{4.0e-2}       & 2.0e-1               & 2.7e-1               & \multicolumn{1}{c|}{2.8e-1}       & 1.0e+0               & 1.1e+3               & \multicolumn{1}{c|}{1.0e+4}       & 1.0e+0               & 1.2e+1               & \multicolumn{1}{c|}{1.0e+2}       \\
                                  &                                             &                      & $\check{f}$                & 5.0e-1               & 5.0e-1               & \multicolumn{1}{c|}{5.0e-1}       & 1.0e+0               & 1.0e+0               & \multicolumn{1}{c|}{1.0e+0}       & 1.0e+0               & 1.0e+0               & \multicolumn{1}{c|}{1.0e+0}       & 1.0e+0               & 1.0e+0               & \multicolumn{1}{c|}{1.0e+0}       \\ \hline
\multirow{2}{*}{\textbf{CCM$_3$}} & \multirow{2}{*}{535}                        & \multirow{2}{*}{4}   & $f$                        & 4.0e-5               & 3.6e-2               & \multicolumn{1}{c|}{4.0e-2}       & 2.0e-1               & 2.7e-1               & \multicolumn{1}{c|}{2.8e-1}       & 1.0e+0               & 1.1e+5               & \multicolumn{1}{c|}{1.0e+6}       & 1.0e+0               & 1.1e+2               & \multicolumn{1}{c|}{1.0e+3}       \\
                                  &                                             &                      & $\check{f}$                & 5.0e-1               & 5.0e-1               & \multicolumn{1}{c|}{5.0e-1}       & 1.0e+0               & 1.0e+0               & \multicolumn{1}{c|}{1.0e+0}       & 1.0e+0               & 1.0e+0               & \multicolumn{1}{c|}{1.0e+0}       & 1.0e+0               & 1.0e+0               & \multicolumn{1}{c|}{1.0e+0}       \\ \hline
\end{tabular}}
	\caption{Properties of faces of meshes used in the Test 3 and in the Test 4.}
	\label{tab:properties_faces_mesh3D}
\end{table}
\begin{table}[]
	\centering
	\resizebox{\textwidth}{!}{
\begin{tabular}{lcccccccccccccclll}
\multirow{2}{*}{}                 & \multirow{2}{*}{\textbf{$\#\mathcal{T}_h$}} & \textbf{$N_{v}^E$}    & \textbf{$N_{e}^E$}    & \textbf{$N_{f}^E$}    & \multirow{2}{*}{\textbf{}} & \multicolumn{3}{c|}{\textbf{Volume}}                                            & \multicolumn{3}{c|}{\textbf{Diameter}}                                          & \multicolumn{3}{c|}{\textbf{Anisotropic  Ratio}}                                & \multicolumn{3}{c|}{\textbf{Face  Ratio}}                                                               \\
                                  &                                             & \textbf{avg}          & \textbf{avg}          & \textbf{avg}          &                            & \textbf{min}         & \textbf{avg}         & \multicolumn{1}{c|}{\textbf{max}} & \textbf{min}         & \textbf{avg}         & \multicolumn{1}{c|}{\textbf{max}} & \textbf{min}         & \textbf{avg}         & \multicolumn{1}{c|}{\textbf{max}} & \multicolumn{1}{c}{\textbf{min}} & \multicolumn{1}{c}{\textbf{avg}} & \multicolumn{1}{c|}{\textbf{max}} \\ \hline
\multirow{2}{*}{\textbf{RTTM}}    & \multirow{2}{*}{569}                        & \multirow{2}{*}{4}    & \multirow{2}{*}{6}    & \multirow{2}{*}{4}    & $E$                        & 3.3e-4               & 1.8e-3               & \multicolumn{1}{c|}{5.0e-3}       & 1.8e-1               & 3.5e-1               & \multicolumn{1}{c|}{4.7e-1}       & 1.5e+0               & 1.0e+1               & \multicolumn{1}{c|}{1.9e+2}       & 1.1e+0                           & 1.8e+0                           & \multicolumn{1}{l|}{3.5e+0}       \\
                                  &                                             &                       &                       &                       & $\hat{E}$                  & 1.2e-1               & 1.2e-1               & \multicolumn{1}{c|}{1.2e-1}       & 1.0e+0               & 1.0e+0               & \multicolumn{1}{c|}{1.0e+0}       & 1.0e+0               & 1.0e+0               & \multicolumn{1}{c|}{1.0e+0}       & 1.0e+0                           & 1.0e+0                           & \multicolumn{1}{l|}{1.0e+0}       \\ \hline
\multirow{2}{*}{\textbf{GPDM}}    & \multirow{2}{*}{308}                        & \multirow{2}{*}{20.5} & \multirow{2}{*}{30.6} & \multirow{2}{*}{12.1} & $E$                        & 6.7e-7               & 3.2e-3               & \multicolumn{1}{c|}{6.2e-3}       & 4.6e-2               & 3.2e-1               & \multicolumn{1}{c|}{7.9e-1}       & 1.4e+0               & 9.2e+0               & \multicolumn{1}{c|}{4.2e+2}       & 1.4e+0                           & 1.8e+3                           & \multicolumn{1}{l|}{8.9e+4}       \\
                                  &                                             &                       &                       &                       & $\hat{E}$                  & 1.2e-1               & 1.8e-1               & \multicolumn{1}{c|}{2.2e-1}       & 1.0e+0               & 1.0e+0               & \multicolumn{1}{c|}{1.0e+0}       & 1.0e+0               & 1.0e+0               & \multicolumn{1}{c|}{1.0e+0}       & 1.0e+0                           & 1.7e+3                           & \multicolumn{1}{l|}{7.5e+4}       \\ \hline
                                  & \multicolumn{1}{l}{}                        & \multicolumn{1}{l}{}  & \multicolumn{1}{l}{}  & \multicolumn{1}{l}{}  & \multicolumn{1}{l}{}       & \multicolumn{1}{l}{} & \multicolumn{1}{l}{} & \multicolumn{1}{l}{}              & \multicolumn{1}{l}{} & \multicolumn{1}{l}{} & \multicolumn{1}{l}{}              & \multicolumn{1}{l}{} & \multicolumn{1}{l}{} & \multicolumn{1}{l}{}              &                                  &                                  &                                   \\ \hline
\multirow{2}{*}{\textbf{CCM$_1$}} & \multirow{2}{*}{150}                        & \multirow{2}{*}{8}    & \multirow{2}{*}{12}   & \multirow{2}{*}{6}    & $E$                        & 8.0e-4               & 6.7e-3               & \multicolumn{1}{c|}{8.0e-3}       & 2.8e-1               & 3.3e-1               & \multicolumn{1}{c|}{3.5e-1}       & 1.0e+0               & 1.8e+1               & \multicolumn{1}{c|}{1.0e+2}       & 1.0e+0                           & 2.5e+0                           & \multicolumn{1}{l|}{1.0e+1}       \\
                                  &                                             &                       &                       &                       & $\hat{E}$                  & 1.9e-1               & 1.9e-1               & \multicolumn{1}{c|}{1.9e-1}       & 1.0e+0               & 1.0e+0               & \multicolumn{1}{c|}{1.0e+0}       & 1.0e+0               & 1.0e+0               & \multicolumn{1}{c|}{1.0e+0}       & 1.0e+0                           & 1.0e+0                           & \multicolumn{1}{l|}{1.0e+0}       \\ \hline
\multirow{2}{*}{\textbf{CCM$_2$}} & \multirow{2}{*}{150}                        & \multirow{2}{*}{8}    & \multirow{2}{*}{12}   & \multirow{2}{*}{6}    & $E$                        & 8.0e-5               & 6.7e-3               & \multicolumn{1}{c|}{8.0e-3}       & 2.8e-1               & 3.4e-1               & \multicolumn{1}{c|}{3.5e-1}       & 1.0e+0               & 1.7e+3               & \multicolumn{1}{c|}{1.0e+4}       & 1.0e+0                           & 1.8e+1                           & \multicolumn{1}{l|}{1.0e+2}       \\
                                  &                                             &                       &                       &                       & $\hat{E}$                  & 1.9e-1               & 1.9e-1               & \multicolumn{1}{c|}{1.9e-1}       & 1.0e+0               & 1.0e+0               & \multicolumn{1}{c|}{1.0e+0}       & 1.0e+0               & 1.0e+0               & \multicolumn{1}{c|}{1.0e+0}       & 1.0e+0                           & 1.0e+0                           & \multicolumn{1}{l|}{1.0e+0}       \\ \hline
\multirow{2}{*}{\textbf{CCM$_3$}} & \multirow{2}{*}{150}                        & \multirow{2}{*}{8}    & \multirow{2}{*}{12}   & \multirow{2}{*}{6}    & $E$                        & 8.0e-6               & 6.7e-3               & \multicolumn{1}{c|}{8.0e-3}       & 2.8e-1               & 3.4e-1               & \multicolumn{1}{c|}{3.5e-1}       & 1.0e+0               & 1.7e+5               & \multicolumn{1}{c|}{1.0e+6}       & 1.0e+0                           & 1.7e+2                           & \multicolumn{1}{l|}{1.0e+3}       \\
                                  &                                             &                       &                       &                       & $\hat{E}$                  & 1.9e-1               & 1.9e-1               & \multicolumn{1}{c|}{1.9e-1}       & 1.0e+0               & 1.0e+0               & \multicolumn{1}{c|}{1.0e+0}       & 1.0e+0               & 1.0e+0               & \multicolumn{1}{c|}{1.0e+0}       & 1.0e+0                           & 1.0e+0                           & \multicolumn{1}{l|}{1.0e+0}       \\ \hline
\end{tabular}}
	\caption{Properties of bulks of meshes used in the Test 3 and in the Test 4.}
	\label{tab:properties_bulk_mesh3D}
\end{table}
\begin{figure}[]
	\centering
	\subfigure[\label{fig:condPi0Derx_Tetra200_Poisson}]
	{\includegraphics[width=.32\textwidth, height = .17\textheight]{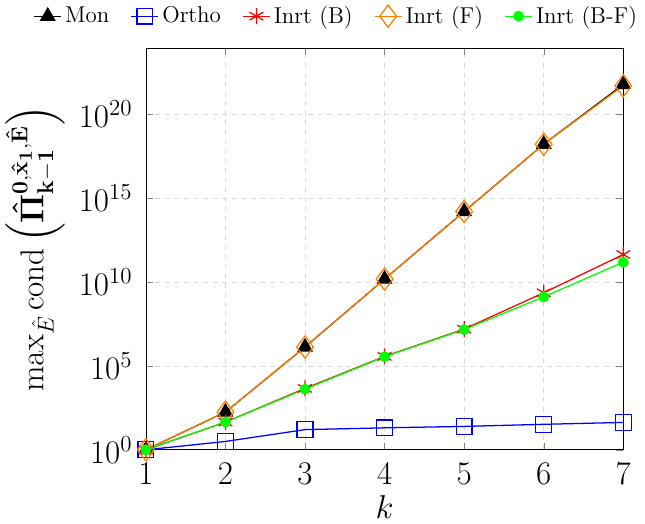}}
	\subfigure[\label{fig:condPiNabla_Tetra200_Poisson}]
	{\includegraphics[width=.32\textwidth, height = .17\textheight]{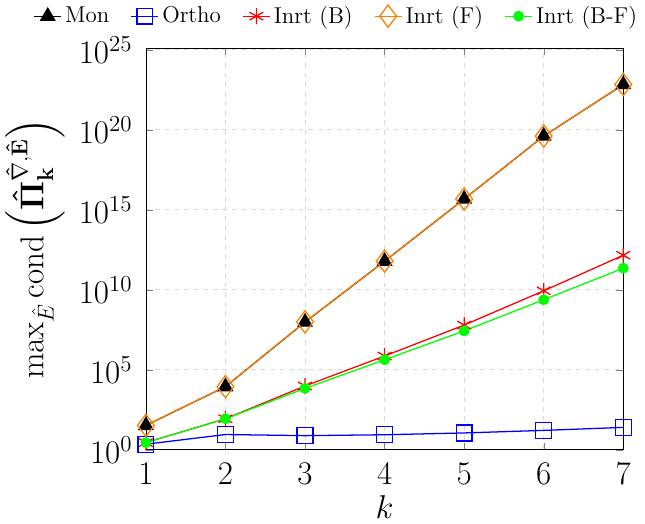}}
	\subfigure[\label{fig:condPi0km1_Tetra200_Poisson}]
	{\includegraphics[width=.32\textwidth, height = .17\textheight]{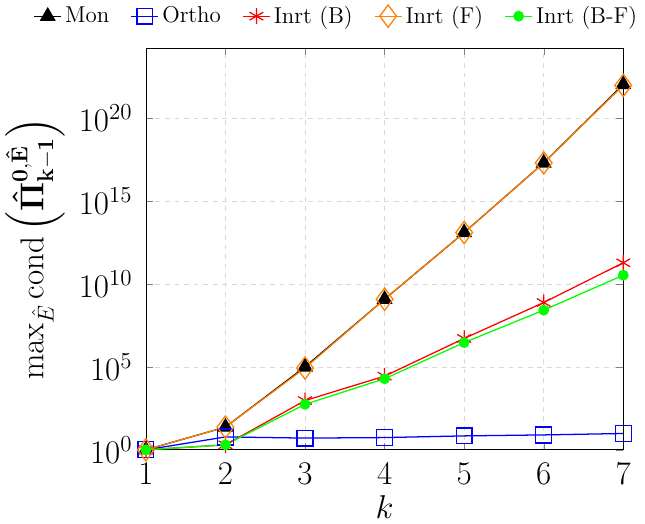}}\\
	\subfigure[\label{fig:condPi0Derx_conf_Poisson}]
	{\includegraphics[width=.32\textwidth, height = .17\textheight]{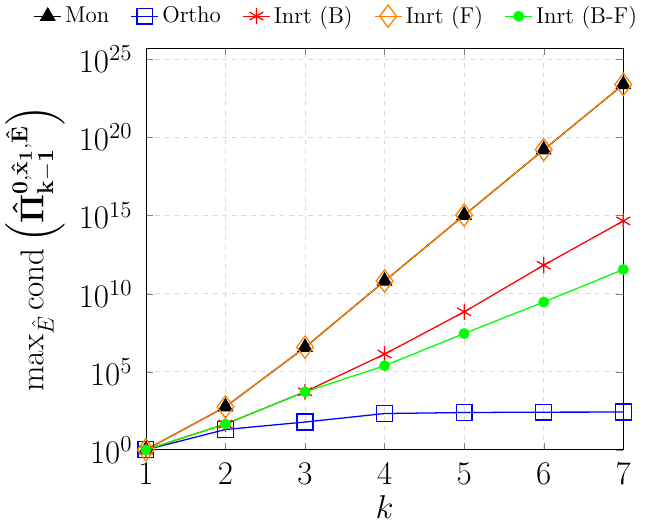}}
	\subfigure[\label{fig:condPiNabla_conf_Poisson}]
	{\includegraphics[width=.32\textwidth, height = .17\textheight]{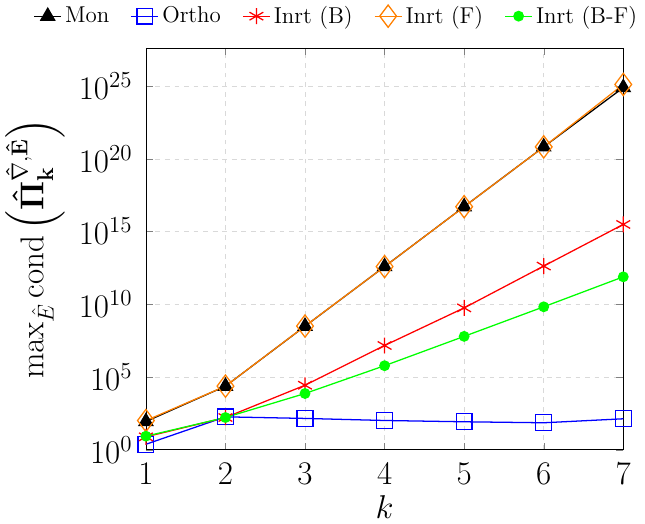}}
	\subfigure[\label{fig:condPi0km1_conf_Poisson}]
	{\includegraphics[width=.32\textwidth, height = .17\textheight]{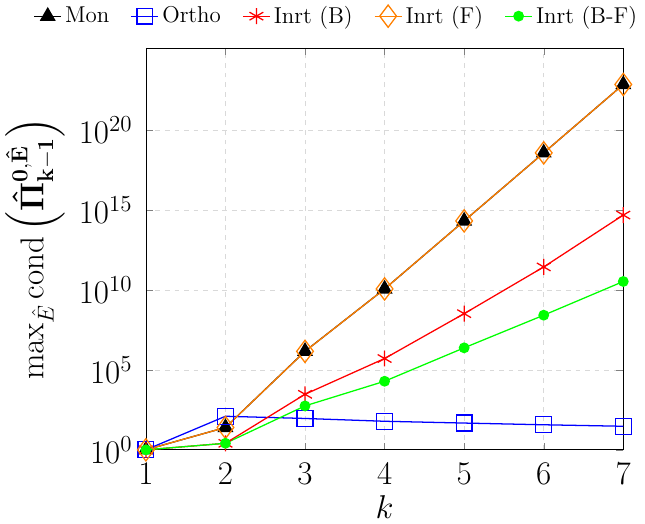}}
	\caption{Test 3: Behaviours of the worst condition numbers of 3D local projection matrices among elements at varying $k$. Top: RTTM. Bottom: GPDM.}
	\label{fig:condlocalprojection_test3}
\end{figure}
\begin{figure}[]
	\centering
	\subfigure[\label{fig:condPiNablaFace_Tetra200_Poisson}]
	{\includegraphics[width=.32\textwidth, height = .17\textheight]{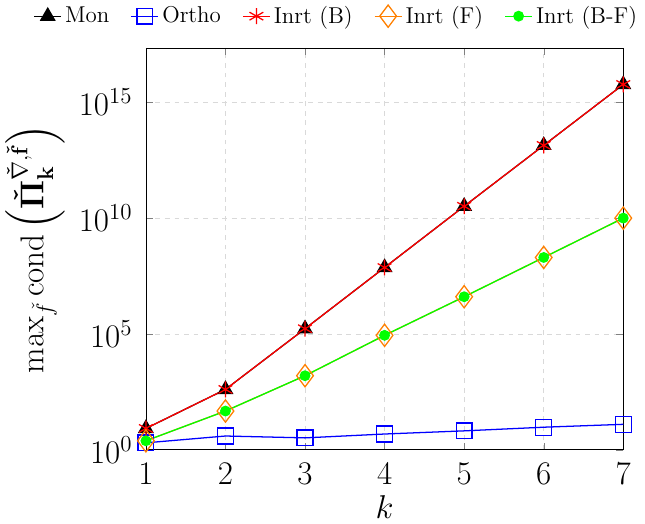}}
	\subfigure[\label{fig:condPi0km1Face_Tetra200_Poisson}]
	{\includegraphics[width=.32\textwidth, height = .17\textheight]{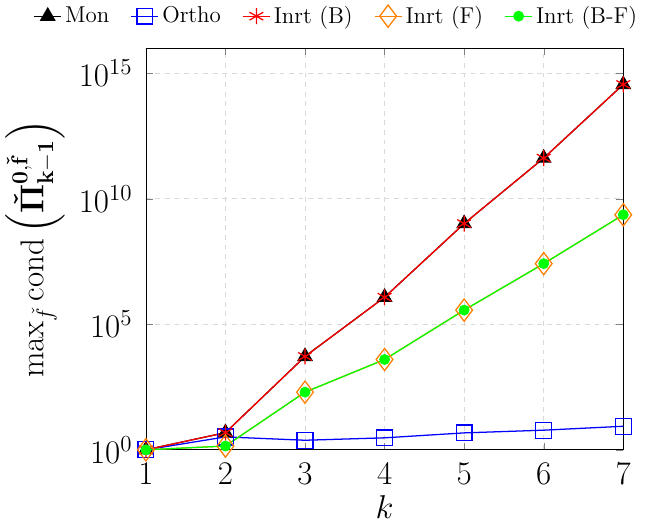}}\\
	\subfigure[\label{fig:condPiNablaFace_conf_Poisson}]
	{\includegraphics[width=.32\textwidth, height = .17\textheight]{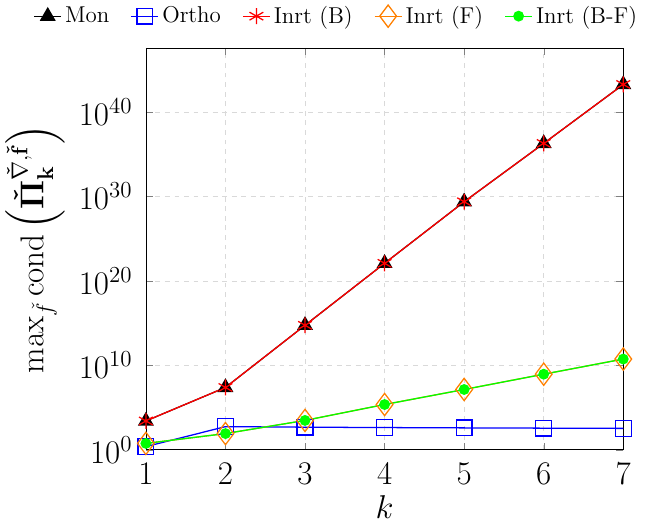}}
	\subfigure[\label{fig:condPi0km1Face_conf_Poisson}]
	{\includegraphics[width=.32\textwidth, height = .17\textheight]{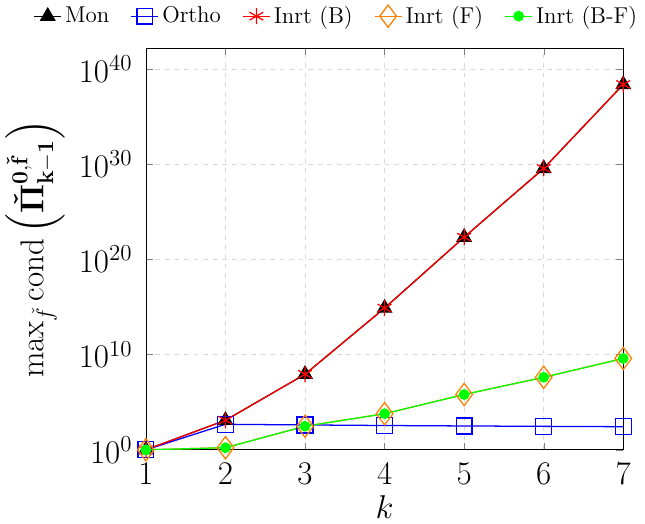}}
	\caption{Test 3: Behaviours of the worst condition numbers of 2D local projection matrices among faces at varying $k$. Top: RTTM. Bottom: GPDM.}
	\label{fig:condlocalprojection_face_test3}
\end{figure}
\begin{figure}[]
	\centering
	\subfigure[\label{fig:condStiff_Tetra200_Poisson}]
	{\includegraphics[width=.32\textwidth, height = .17\textheight]{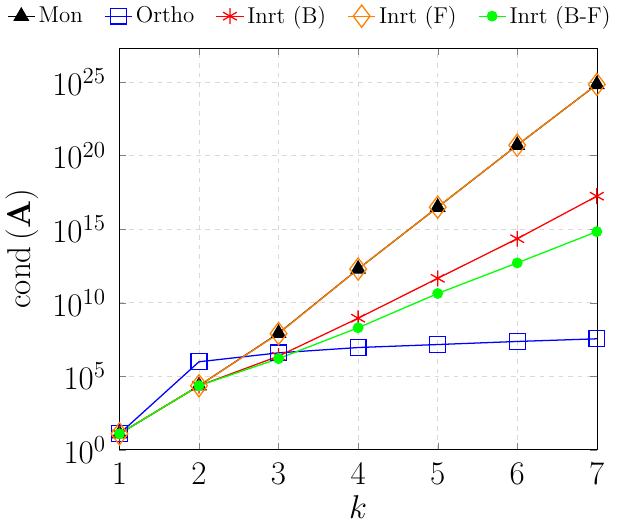}}
	\subfigure[\label{fig:ErrorL2_Tetra200_Poisson}]
	{\includegraphics[width=.32\textwidth, height = .17\textheight]{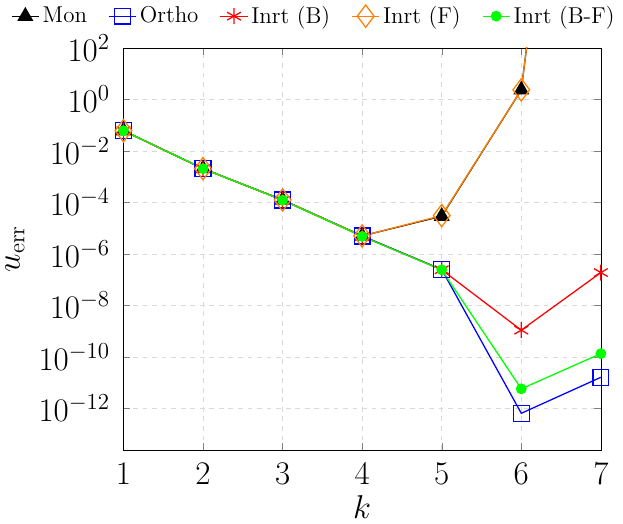}}
	\subfigure[\label{fig:ErrorH1_Tetra200_Poisson}]
	{\includegraphics[width=.32\textwidth, height = .17\textheight]{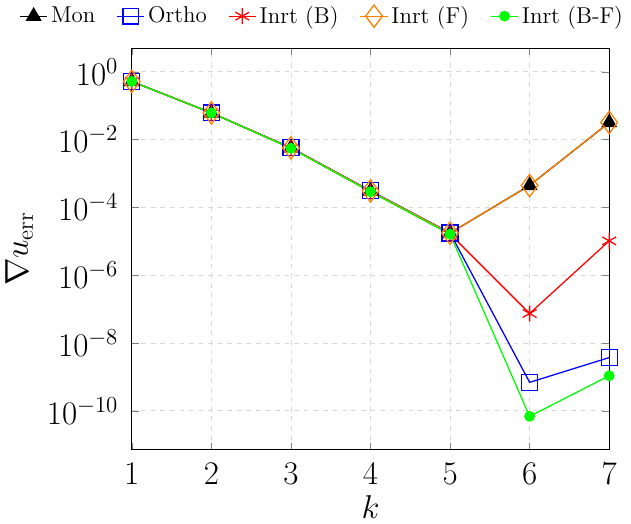}}\\
	\subfigure[\label{fig:condStiff_conf_Poisson}]
	{\includegraphics[width=.32\textwidth, height = .17\textheight]{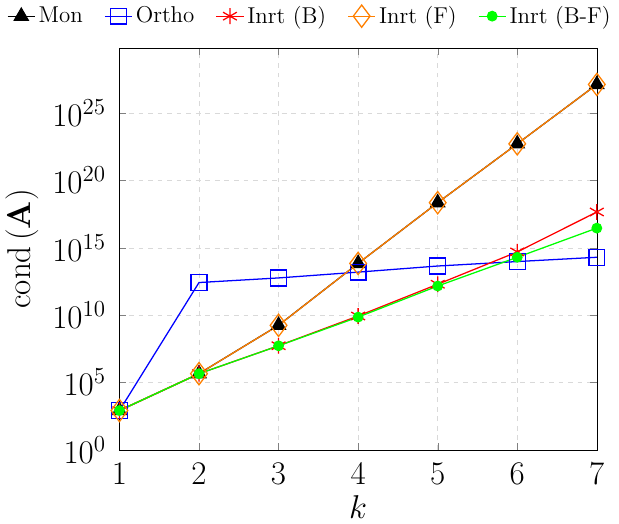}}
	\subfigure[\label{fig:ErrorL2_conf_Poisson}]
	{\includegraphics[width=.32\textwidth, height = .17\textheight]{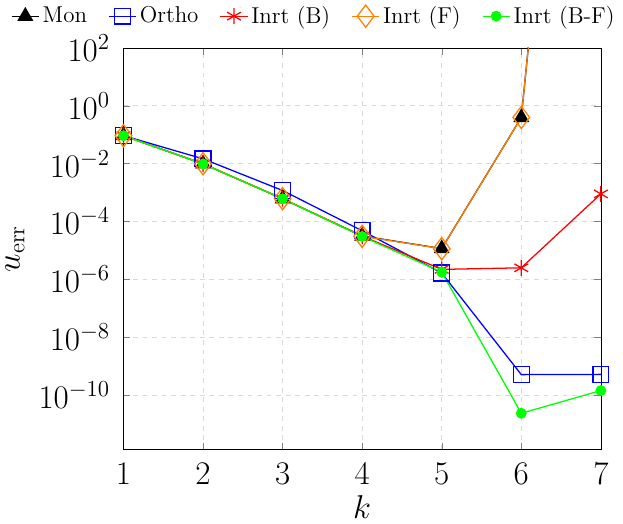}}
	\subfigure[\label{fig:ErrorH1_conf_Poisson}]
	{\includegraphics[width=.32\textwidth, height = .17\textheight]{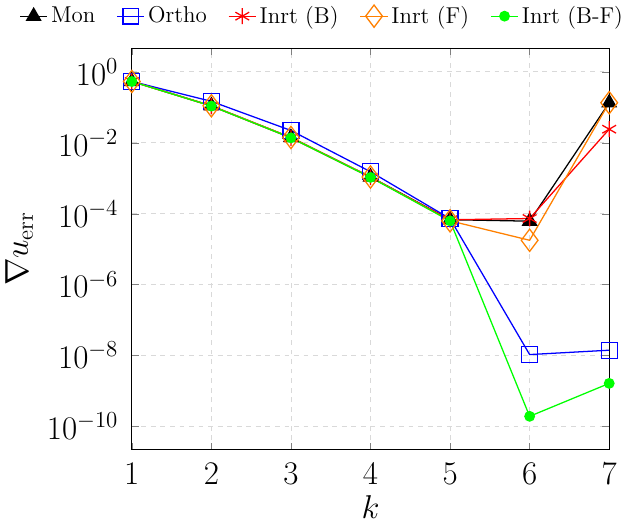}}
	\caption{Test 3: Behaviours of the condition number of the global system matrix and of errors \eqref{eq:errorL2} and \eqref{eq:errorH1} at varying $k$. Top: RTTM. Bottom: GPDM.}
	\label{fig:error_cond_test3}
\end{figure}

In the following, we analyze the three-dimensional case.
Let $\Omega =\left(0,1\right)^3$, we consider the problem \eqref{eq:modelproblem} with constant coefficients $\D = \mathbf{I}$, $\gamma = 0$ and $\bb= \bm{0}$, where we defined $f$ and the non-homogeneous Dirichlet boundary condition in such a way the exact solution is the polynomial function of degree $6$:

\begin{equation*}
	u(x_1,x_2,x_3) = 1.7 + 64 x_1 x_2 x_3(1-x_1)(1-x_2)(1-x_3).
\end{equation*}

In this experiment, we generate:
\begin{itemize}
	\item a regular tetrahedral mesh (RTTM in short) to validate our procedures in the 3D case (Figure~\ref{fig:mesh_Tetra200});
	\item a generic polyhedral mesh (GPDM in short) made up mostly of badly-shaped polyhedrons and characterized by a copious number of aligned faces (Figure~\ref{fig:mesh_conf}).
\end{itemize}

The main geometric properties of the polyhedrons and of the polygonal faces belonging to these meshes are shown in Tables~\ref{tab:properties_faces_mesh3D} and \ref{tab:properties_bulk_mesh3D}. More precisely, for the polyhedron elements we consider the volume, the diameter, the anisotropic ratio and the \textit{face ratio}, i.e. the ratio between the highest and the smallest areas of the faces of a polyhedron; for the polygonal faces we measure the area, the diameter, the anisotropic ratio and the edge ratio.
As in the 2D case, the application of $F_E$ map tends to uniform elements belonging to the same categories, as happens in the case of tetrahedrons. Moreover, we note that the face ratio of the mapped elements is not far from the face ratio of the original elements, since the map does not take into account the presence of aligned faces. However, we do not care about the geometric properties related to the faces of the mapped elements $\srescale{E}$, since we never them (see Remark \ref{rem:notmappedfaces}).

We measure the performances of the different VEM approaches up to polynomial degree $k=7$, as the dimension of the local and global matrices increase faster in the three-dimensional setting. 

The local condition numbers are shown in Figure~\ref{fig:condlocalprojection_test3} for the 3D projectors and in Figure~\ref{fig:condlocalprojection_face_test3} for the 2D projectors on faces. We further omit the graphs reporting the behaviours of the condition numbers of $\mmproj{0,\srescale{x}_2,\srescale{E}}{k-1}$ and $\mmproj{0,\srescale{x}_3,\srescale{E}}{k-1}$ at varying $k$, since these trends are very similar to the one of the condition number of $\mmproj{0,\srescale{x}_1,\srescale{E}}{k-1}$ for each method.

By looking at those figures, we can observe that by using the inertial approach \ref{itm:F}, and thus by mapping only the faces of the elements, we improve just the condition numbers of the 2D local projection matrices, whereas the condition numbers of the 3D local projection matrices are comparable to the ones obtained with the standard monomial approach. For the same reason, the condition numbers of 2D local projection matrices related to the inertial approach \ref{itm:B} are in the same order of magnitude as the monomial ones. We appreciate that the condition numbers of the 3D local projection matrices in the inertial approach \ref{itm:B} differ from the ones related to the inertial approach \ref{itm:BF} even if we use the same mapping for their bulks, because of the different contributions of the boundary integrals. Finally, we note that the condition numbers of 2D local projection matrices are equal to the ones obtained by resorting to the inertial approach \ref{itm:F}, since the projectors are computed on the same faces in the two cases.

The Figure~\ref{fig:error_cond_test3} shows the trends of the condition number of the system matrix and of the errors \eqref{eq:errorL2} and \eqref{eq:errorH1} at varying $k$.
As expected from the local behaviour, the performances of the inertial approach \ref{itm:F} overlap the monomial one in both tests.
On the other hand, the other two inertial approaches \ref{itm:B} and \ref{itm:BF} lead to very good results in terms of the condition number of the global matrix.
Moreover, we observe that the global performances of the \ref{itm:BF} strategy are comparable and sometimes better than the ones obtained with the orthonormal approach in the case of polyhedral mesh GPDM.
Finally, the errors of \ref{itm:BF} and of orthonormal methods correctly decay to zero when $k=6$, since the exact solution is a polynomial function of degree $6$, whereas the errors related to the other approaches start to raise for $k\geq 5$.

In conclusion, the results show that mapping only the faces or only the bulks of elements of the mesh is not sufficient to obtain more reliable and accurate solutions. On the other hand, using the inertial approach \ref{itm:BF} is very efficient from a computational point of view and leads to very high-quality results.

\subsection{Test 4: Collapsing polyhedrons}
\begin{figure}[]
	\centering
	\subfigure[\label{fig:condPi0Derx_Hexahedral_Elliptic}]
	{\includegraphics[width=.32\textwidth, height = .17\textheight]{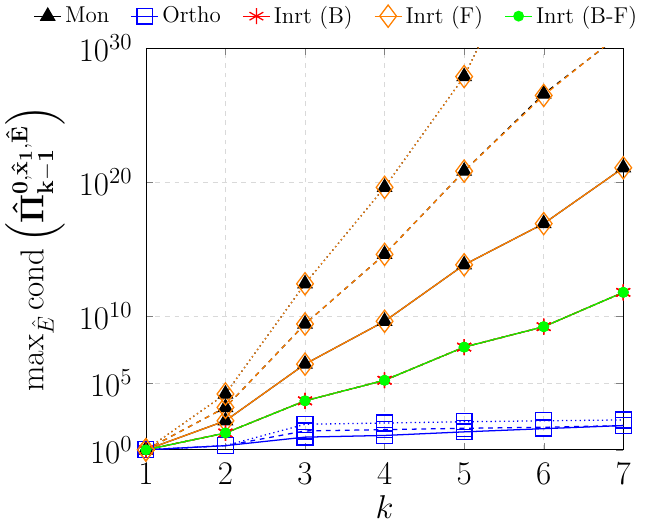}}
	\subfigure[\label{fig:condPi0Dery_Hexahedral_Elliptic}]
	{\includegraphics[width=.32\textwidth, height = .17\textheight]{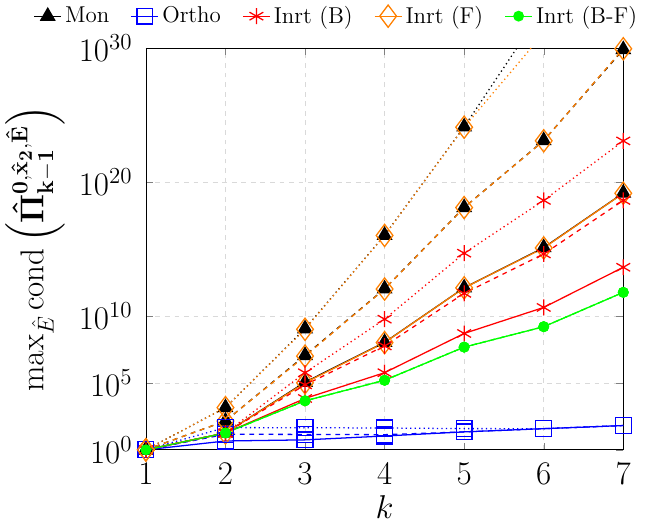}}
	\subfigure[\label{fig:condPi0Derz_Hexahedral_Elliptic}]
	{\includegraphics[width=.32\textwidth, height = .17\textheight]{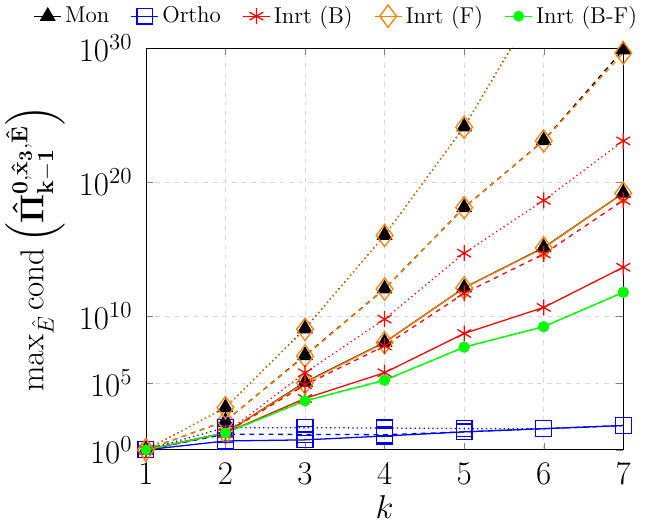}}
	\subfigure[\label{fig:condPiNabla_Hexahedral_Elliptic}]
	{\includegraphics[width=.32\textwidth, height = .17\textheight]{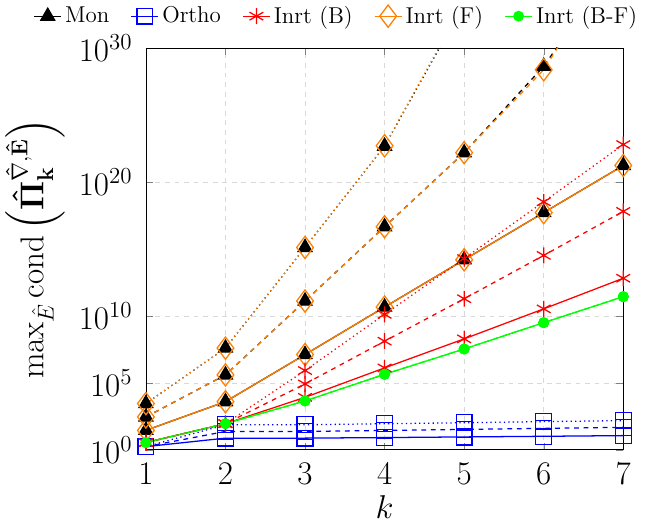}}
	\subfigure[\label{fig:condPi0km1_Hexahedral_Elliptic}]
	{\includegraphics[width=.32\textwidth, height = .17\textheight]{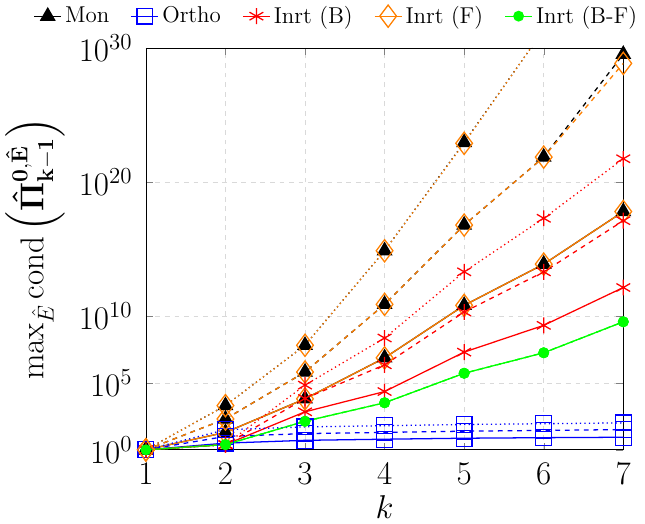}}
	\caption{Test 4: Behaviours of the worst condition numbers of 3D local projection matrices among elements at varying $k$. Solid lines: CCM$_1$. Dashed lines: CCM$_2$. Dotted lines: CCM$_3$.}
	\label{fig:condlocalprojection_test4}
\end{figure}
\begin{figure}[]
	\centering
	\subfigure[\label{fig:condPiNablaFace_Hexahedral_Elliptic}]
	{\includegraphics[width=.32\textwidth, height = .17\textheight]{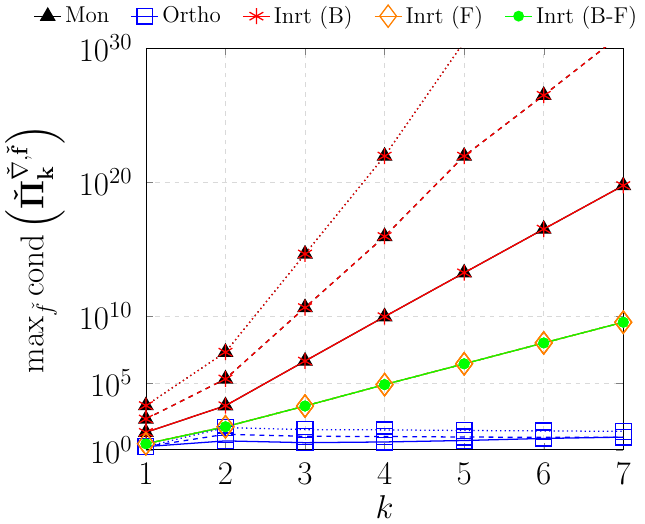}}
	\subfigure[\label{fig:condPi0km1Face_Hexahedral_Elliptic}]
	{\includegraphics[width=.32\textwidth, height = .17\textheight]{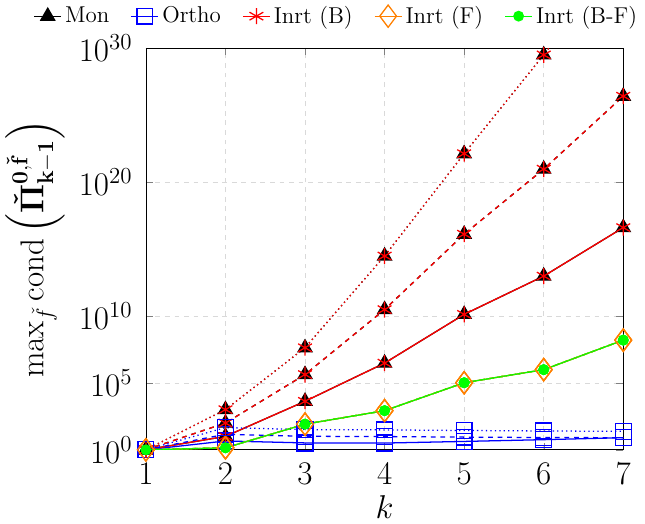}}
	\caption{Test 4: Behaviours of the worst condition numbers of 2D local projection matrices among faces at varying $k$. Solid lines: CCM$_1$. Dashed lines: CCM$_2$. Dotted lines: CCM$_3$.}
	\label{fig:condlocalprojection_face_test4}
\end{figure}
\begin{figure}[]
	\centering
	\subfigure[\label{fig:ErrorL2_Hexahedral_Elliptic}]
	{\includegraphics[width=.32\textwidth, height = .17\textheight]{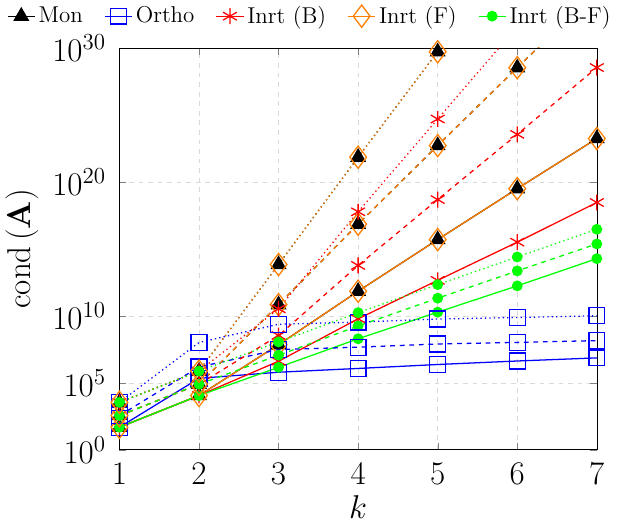}}
	\subfigure[\label{fig:ErrorL2_Hexahedral_Elliptic}]
	{\includegraphics[width=.32\textwidth, height = .17\textheight]{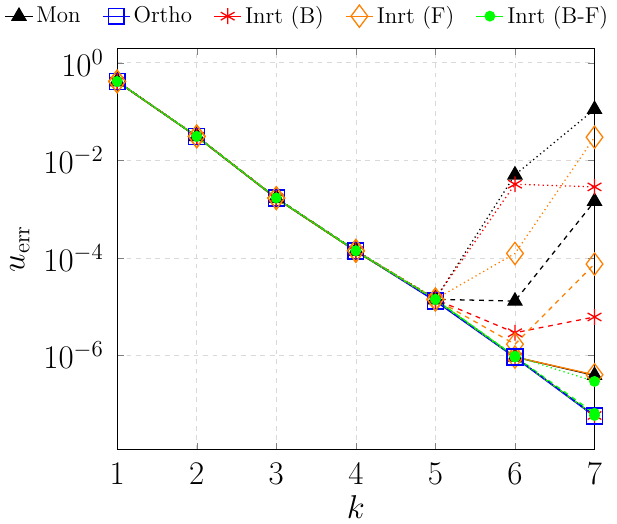}}
	\subfigure[\label{fig:ErrorH1_Hexahedral_Elliptic}]
	{\includegraphics[width=.32\textwidth, height = .17\textheight]{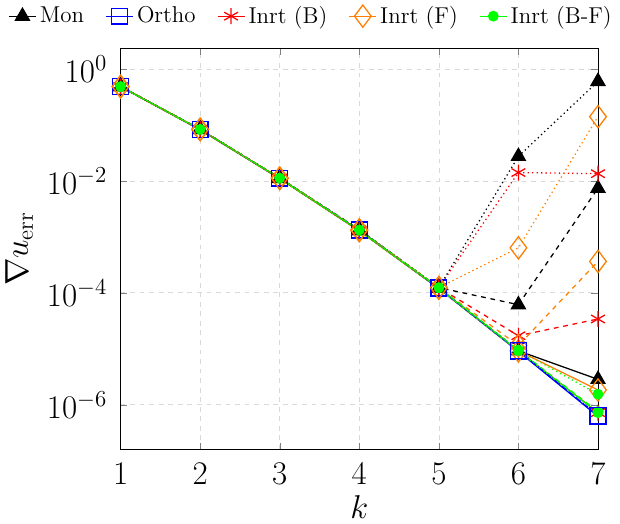}}
	\caption{Test 3: Behaviours of the condition number of the global system matrix and of errors \eqref{eq:errorL2} and \eqref{eq:errorH1} at varying $k$. Solid lines: CCM$_1$. Dashed lines: CCM$_2$. Dotted lines: CCM$_3$.}
	\label{fig:error_cond_test4}
\end{figure}

This last experiment is a natural extension of Test 2 to the 3D case.
Thus, we consider on $\Omega =\left(0,1\right)^3$ the advection-diffusion-reaction problem \eqref{eq:modelproblem} with variable coefficients given by

\begin{equation*}
	\D(x_1,x_2,x_3) = \begin{bmatrix}
		1 + x_2^2 + x_3^2 & -x_1x_2         & - x_1x_3        \\
		-x_1x_2           & 1+x_1^2 + x_3^2 & - x_2x_3        \\
		-x_1x_3           & -x_2x_3         & 1 + x_1^2+x_2^2
	\end{bmatrix},
\end{equation*}
\begin{equation*}
	\bb(x_1,x_2,x_3) = \begin{bmatrix}
		x_1 \\ x_2 \\ -2x_3
	\end{bmatrix}, \quad \gamma(x_1,x_2,x_3) = x_1x_2x_3.
\end{equation*}
and we set $f$ in such a way the exact solution is

\begin{equation*}
	u(x_1,x_2,x_3) = \sin(\pi x_1) \sin(\pi x_2)  \sin(\pi x_3).
\end{equation*}

We further generate a sequence of three hexahedral meshes with cubic elements of edge $0.2$ with the exception of a central band along $x_1$-axis made by two groups of hexahedrons, one of which (the purple band highlighted in Figure~\ref{fig:mesh_Hexahedron_1e_01}) is composed by hexahedrons of volume which varies from $8 \cdot 10^{-4}$ in the first mesh CCM$_1$ to $8 \cdot 10^{-6}$ in the last mesh of the sequence CCM$_3$.
Figure~\ref{fig:mesh_Hexahedron_1e_01} shows a clip of the first mesh, and Tables~\ref{tab:properties_faces_mesh3D} and \ref{tab:properties_bulk_mesh3D} report the geometrical 2D and 3D properties of the mesh on rows CCM$_i$, $i \in \{1,2,3\}$. 

Figures ~\ref{fig:condlocalprojection_test4} and \ref{fig:condlocalprojection_face_test4} report the condition numbers of local projection matrices in the three meshes at varying polynomial degree $k$.
For this particular test, we decide to report on the plots the behaviour of the worst local condition number of all the three directions of projector $\mproj{0,\hat{x}_i,\srescale{E}}{k-1}$ to highlight the differences in the $x_1$-axis direction of the central band with respect the other directions.


As happens for the 2D case, mapping only the faces with the \ref{itm:F} approach makes uniform the condition numbers of the 2D local projection matrices among meshes in the sequence. On the other hand, mapping only the bulks in the \ref{itm:B} approach does not make the condition numbers of the 3D local projection matrices independent of the mesh of the sequence because the boundary contributions depend on the actual mesh.
Finally, mapping both faces and bulks with the inertial \ref{itm:BF} approach makes the condition numbers of all local projection matrices approximately independent of the features of the central band.

To conclude the analysis, we show the global performances in terms of the condition number of the global matrix and of the errors \eqref{eq:errorL2} and \eqref{eq:errorH1} in Figure~\ref{fig:error_cond_test4}. As concluded for the Test $3$, the global performances of all the inertial approaches depend on the geometric properties of the original meshes of the sequence.
Furthermore, the inertial approach \ref{itm:BF} reveals to be more robust and more accurate than the monomial approach and it has a behaviour comparable to the one of the orthonormal approach, besides being the best advantageous choice from a computational point of view.

\section{Conclusions}
One of the main features of the Virtual Element Method is the possibility to use very complex geometries, but the use of the classical scaled monomial basis in its construction generally leads to very low-quality results in presence of badly-shaped polytopes and for high polynomial degrees.

In this paper, we propose a new procedure to build a polynomial basis on polytopes to mitigate the ill-conditioning of the local projection matrices and of the global system matrix with a computational complexity which depends only on the geometric dimension of the problem ($d=2$ or $3$) and on the number of faces in 3D.

Throughout different numerical experiments of increasing complexity in 2D and 3D cases, we observe that recomputing the scaled monomial basis on more well-shaped polytopes leads to obtain an acceptable well-conditioned polynomial basis and, consequently, a more reliable solution.
We highlight the need of defining both the 2D and 3D local polynomial bases on well-shaped polytopes in the 3D case in order to improve the local and global performances with respect to the use of the standard scaled monomial basis.

Finally, the proposed approach has proved to have reasonable and, most of the time, comparable results from a practical point of view with respect to the use of a local orthonormal polynomial basis, while being less expensive from a computational point of view.
\section*{Acknowledgments}

This work is supported by INdAM-GNCS, the MIUR project “Dipartimenti di Eccellenza 2018-2022” (CUP E11G18000350001), by PRIN project “Advanced polyhedral discretisations of heterogeneous PDEs for multiphysics problems” (0204LN5N5\_003) and by the MIUR programme ``Programma Operativo Nazionale Ricerca e Innovazione 2014-2020'' (D.M. 1061/2021).
Computational resources are partially supported by SmartData@polito.

\bibliographystyle{IEEEtran}
\bibliography{biblio.bib}

\end{document}